\documentclass[12pt,letterpaper]{amsart}
\usepackage[margin=1in]{geometry}

\usepackage{amsmath,amssymb,amsthm,mathrsfs}
\usepackage{enumitem,booktabs,url}
\usepackage{xcolor}
\definecolor{referenceblue}{RGB}{0,56,105}
\usepackage[colorlinks=true,linkcolor=referenceblue,citecolor=referenceblue,
            urlcolor=referenceblue,filecolor=referenceblue]{hyperref}
\hypersetup{pdftitle={Direct virial estimates and asymptotic stability of gray solitons},
pdfauthor={Justin Holmer and Numann Malik}}
\newcommand{\sech}{\operatorname{sech}}
\newcommand{\spn}{\operatorname{span}}

\renewcommand{\Re}{\operatorname{Re}}
\renewcommand{\Im}{\operatorname{Im}}
\newtheorem{theorem}{Theorem}[section]
\newtheorem{proposition}[theorem]{Proposition}
\newtheorem{lemma}[theorem]{Lemma}
\newtheorem{corollary}[theorem]{Corollary}

\theoremstyle{remark}

\allowdisplaybreaks[1]
\numberwithin{equation}{section}
\title{Direct virial estimates and asymptotic stability of gray solitons}
\author{Justin Holmer}
\address{Department of Mathematics, Brown University,
Providence, RI 02912}
\email{justin\_holmer@brown.edu}
\author{Numann Malik}
\address{Department of Mathematics and Statistics, The College of New Jersey,
Ewing, NJ 08628}
\email{malikn@tcnj.edu}
\date{\today}
\begin{document}
\begin{abstract}
We study gray solitons of the one-dimensional defocusing Gross--Pitaevskii
equation. Their asymptotic stability in the energy space is known from
the work of B\'ethuel, Gravejat and Smets. Here we obtain a direct
quantitative estimate on the evolving perturbation by coupling a
localized momentum identity with a regularized dual virial identity.
For the dual identity we prove constrained coercivity with the weight
$y(1-\hat\eta_\beta(y))$ by a half-line argument and an exact polynomial
certificate. Spatial truncation and frequency regularization yield
a time-integrable local energy, finite total variation of the fitted
soliton angle, and square-integrable error in the center velocity.
For the unforced cubic equation, the conserved energies of Koch and
Liao give uniform control of the high frequencies along each solution.
It follows that every small energy-space perturbation converges to zero
strongly in $H^1_{\mathrm{loc}}\times L^2_{\mathrm{loc}}$ in the fitted
frame and that the velocity of this frame tends to the limiting soliton
speed.
\end{abstract}
\maketitle

\section{Introduction}
We consider the defocusing Gross--Pitaevskii equation
\begin{equation}\label{E:GP}
i u_t+u_{xx}+2(1-|u|^2)u=0,
\qquad |u(x,t)|\longrightarrow1\quad (|x|\longrightarrow\infty).
\end{equation}
The classical dark-soliton family studied by Tsuzuki \cite{Tsuzuki1971}
has the form
\begin{equation}\label{E:q-def}
u(x,t)=q_\beta(x-a_0-2t\sin\beta),\qquad
q_\beta(y)=\cos\beta\tanh(y\cos\beta)+i\sin\beta,
\quad |\beta|<\pi/2.
\end{equation}
The stationary profile $q_0=\tanh y$ is the black soliton. The nonzero
parameters give gray solitons. We work on a fixed gray branch, away
from both the black and sonic endpoints.

Our motivation comes from the effective dynamics of a gray soliton in a
slowly varying momentum field, studied in the companion paper
\cite{HMEffectiveDynamics}. The value of the field at the center
transports the soliton, and its gradient changes the soliton momentum.
Tracking the center and depth on long finite time intervals requires
an estimate for the radiation near the core while the parameters are
moving. This is why we keep both localized virial identities on the
evolving perturbation. The companion uses their spatial estimates in
the forced modulation and energy argument.

For a nonvanishing solution, write $u=|u|e^{i\phi}$ and introduce
\begin{equation}\label{E:hydro}
\eta=1-|u|^2,\qquad
\nu=-\phi_x=-\frac{\Im(\bar u u_x)}{|u|^2}.
\end{equation}
The corresponding soliton, centered at the origin, is
\begin{equation}\label{E:hydro-profile}
\begin{bmatrix}\hat\eta_\beta(y)\\\hat\nu_\beta(y)\end{bmatrix}
=\begin{bmatrix}
\cos^2\beta\,\sech^2(y\cos\beta)\\[2pt]
\displaystyle\frac{\sin\beta\,\cos^2\beta\,\sech^2(y\cos\beta)}
{1-\cos^2\beta\,\sech^2(y\cos\beta)}
\end{bmatrix}.
\end{equation}
We use the energy space $\mathcal X=H^1(\mathbb R)\times L^2(\mathbb R)$,
with $\|F\|_{\mathcal X}^2=\|F_1\|_{H^1}^2+\|F_2\|_2^2$.
In a moving frame $y=x-a(t)$ the perturbation is always denoted by
\begin{equation}\label{E:decomposition}
\begin{bmatrix}\hat\theta(y,t)\\\hat\xi(y,t)\end{bmatrix}
=\begin{bmatrix}\eta(y+a(t),t)\\\nu(y+a(t),t)\end{bmatrix}
-\begin{bmatrix}\hat\eta_{\beta(t)}(y)\\\hat\nu_{\beta(t)}(y)\end{bmatrix}.
\end{equation}
Thus hats on the perturbation indicate the moving coordinates; they
do not denote a dual transformation or a Fourier transform.

\begin{theorem}[Direct estimates and asymptotic stability]\label{T:main}
Fix $0<\beta_0<\pi/2$ and a compact interval
$I\Subset(0,\pi/2)$ with $\beta_0$ in its interior. There exist
$\delta_0>0$, $A\geq1$, and $C>0$ such that the following holds.
Suppose $0<\delta\leq\delta_0$ and
\[
\left\|\begin{bmatrix}
\eta(\cdot,0)-\hat\eta_{\beta_0}(\cdot-a_0)\\
\nu(\cdot,0)-\hat\nu_{\beta_0}(\cdot-a_0)
\end{bmatrix}\right\|_{\mathcal X}\leq\delta.
\]
The solution is global and nonvanishing. There are $C^1$ functions
$\beta(t)\in I$ and $a(t)$, determined by the orthogonality conditions
\eqref{E:orth}, for which
\begin{equation}\label{E:main-estimate}
\begin{aligned}
&\sup_{t\geq0}
\left\|\begin{bmatrix}\hat\theta(t)\\\hat\xi(t)\end{bmatrix}
\right\|_{\mathcal X}^2
+\int_0^\infty\!\int_{\mathbb R}\sech^2(y/A)
\bigl(\hat\theta_y^2+\hat\theta^2+\hat\xi^2\bigr)\,dy\,dt
\leq C\delta^2,\\
&\int_0^\infty|\dot\beta(t)|\,dt
+\int_0^\infty|\dot a(t)-2\sin\beta(t)|^2\,dt\leq C\delta^2.
\end{aligned}
\end{equation}
There is $\beta_\infty\in I$ such that $\beta(t)\to\beta_\infty$ and
$|\beta_\infty-\beta(0)|\le C\delta^2$. For every $R>0$, as $t\to\infty$,
\begin{equation}\label{E:main-strong}
\|\hat\theta(t)\|_{H^1(-R,R)}+\|\hat\xi(t)\|_{L^2(-R,R)}\longrightarrow0,
\qquad \dot a(t)\longrightarrow2\sin\beta_\infty .
\end{equation}
In particular, in the frame selected by \eqref{E:orth} the hydrodynamic
variables converge strongly on bounded intervals to the profile
$(\hat\eta_{\beta_\infty},\hat\nu_{\beta_\infty})$.
The global orbital bound also gives
$(\hat\theta(t),\hat\xi(t))\rightharpoonup(0,0)$ in $\mathcal X$.
The analogous assertions hold on the negative gray branch by the
spatial reflection $(\eta,\nu)(x)\mapsto(\eta(-x),-\nu(-x))$.
\end{theorem}

The constants are locally uniform on a gray branch, subject to keeping
the reference parameter a fixed distance from the boundary of $I$.
They are not asserted to be uniform as $\beta\to0$ or
$|\beta|\to\pi/2$. The bounds in \eqref{E:main-estimate} are quantitative
in the initial perturbation size and independent of the time interval.
The compactness argument for \eqref{E:main-strong} uses the
high-frequency tail of the individual initial datum. It gives a limit
without a uniform pointwise decay rate. The speed conclusion allows a
sublinear drift in $a(t)-2t\sin\beta_\infty$.

\subsection{Orbital and asymptotic stability: historical context}
The stability theory has developed along variational, spectral, and
inverse-scattering approaches. We describe the results for the
one-dimensional scalar cubic equation, and then indicate closely
related extensions. In some of the literature the word ``dark'' is
reserved for the nonzero-speed waves; our term ``gray'' removes that
ambiguity.

The inverse-scattering description of dark solitons and their
interactions goes back to Zakharov and Shabat \cite{ZakharovShabat1973}.
For general nonlinearities, Kivshar and Kr\'olikowski
\cite{KivsharKrolikowski1995} studied a momentum-slope criterion for instability,
and Barashenkov \cite{Barashenkov1996} established the slope criterion
for linear stability.

Early work of Zhidkov on nonzero boundary conditions and kink stability
\cite{Zhidkov1992,Zhidkov2001} established part of the functional
framework for this problem. For gray solitons, Lin \cite{Lin2002}
proved orbital stability using the nonvanishing hydrodynamic variables
and the energy--momentum criterion. Gallo \cite{Gallo2004} justified the
Cauchy theory and conservation laws used in this approach. B\'ethuel, Gravejat and Saut
\cite[Section~2.4]{BGSaut2008} gave a variational proof for the cubic
equation using compactness of energy minimizers at fixed momentum.
The black soliton requires a
different treatment because its amplitude vanishes. Di Menza and
Gallo \cite{DiMenzaGallo2007} proved its linear stability in the cubic
case, within a more general criterion that also detects instability
for other nonlinearities. The full nonlinear orbital stability of the
black soliton under energy-space perturbations was proved by
B\'ethuel, Gravejat, Saut and Smets \cite{BGS_Orbital_2008}.
G\'erard and Zhang \cite{GerardZhang2009} independently obtained orbital
stability near the black soliton by inverse scattering for smooth,
spatially decaying perturbations, including uniform control in the
supremum norm.

Chiron \cite{Chiron2013} developed energy-space stability and instability
criteria for general one-dimensional defocusing nonlinearities,
including the kink endpoint and degenerate energy--momentum diagrams.
B\'ethuel, Gravejat and Smets \cite{BGSChains2014} proved orbital
stability of well-separated chains of gray solitons with distinct
nonzero speeds, ordered according to their speeds. Gallay and
Pelinovsky \cite{GP2015} gave a further proof of black-soliton orbital
stability using a higher conserved quantity and an unconstrained
variational characterization, with perturbations in a topology
involving two derivatives. More recently, Duan and Wang
\cite{DuanWang2025,DuanWangCorrection2025} gave a unified variational
proof for black and gray solitons in the energy framework.

The asymptotic stability of gray solitons under small energy-space
perturbations was proved by B\'ethuel, Gravejat and Smets
\cite{BGS_Asymp_Dark_2015}. Their Theorem~4 gives weak convergence in $H^1\times L^2$ in a frame
whose velocity tends to the limiting soliton speed. Their Theorem~2
also gives local uniform convergence of the complex field after a phase
correction whose derivative tends to zero. Immediately after
Theorem~2 they ask whether the weak convergence of the phase-corrected
gradient can be strengthened to strong convergence on bounded
intervals. Corollary~\ref{C:BGS-gradient} obtains this conclusion from
their theorem and the conserved energies of Koch and Liao
\cite{KochLiaoLow}. The B\'ethuel--Gravejat--Smets proof combines
modulation, localized momentum monotonicity, compactness of limiting
solutions, and a Liouville theorem.
Gravejat and Smets \cite{GS_Asymp_Black_2015}
introduced coordinates and a stability framework that also allow
zeros, and proved asymptotic stability of the black soliton.
In the integrable setting, Vartanian \cite{Vartanian2002} derived
long-time expansions along dark-soliton trajectories under decay and
scattering-data hypotheses. Mohamad \cite{MohamadMarchenko2015} analyzed
the Marchenko system and the large-time behavior of the exact
multisoliton family. Cuccagna and Jenkins
\cite[Theorems~1.1 and~1.6]{CuccagnaJenkins2016} obtained a uniform
description in closed subsectors of the solitonic cone and a stability
theorem for compatible multisolitons under weighted, regular
perturbations. Their opposite far-field values include the single
black soliton, and Remark~1.8 discusses other boundary phases.
Liu and Xu \cite{LiuXu2026} studied partial-mass asymptotics under
weaker weighted hypotheses.
Kowalczyk and Martel \cite{KowalczykMartel2025} subsequently gave a new
proof of the Liouville theorem for both black and gray solitons,
using a factorization of the linearized operator that simplifies the
spectral analysis. They also point toward a direct asymptotic-stability
argument based on that transformation.

For general nonlinearities, Berthoumieu established orbital stability of
ordered, well-separated chains near the sound speed
\cite{BerthoumieuOrbital2025} and proved asymptotic stability for a single
wave in this regime \cite{BerthoumieuSingle2025}. His subsequent preprint
\cite{BerthoumieuChains2026} treats asymptotic stability of well-prepared
chains and constructs solutions asymptotic to prescribed chains.

Thus the stability assertion itself is established literature. The
purpose here is to obtain the estimates \eqref{E:main-estimate} on the
original modulated solution, using the particular weight
$y\zeta_\beta(y)$ defined below, while retaining control of the
regularity, the cutoffs, and the parameter errors. Results for other
nonlinearities or external potentials address additional questions;
for example, the quintic black soliton is treated by Alejo and Corcho
\cite{AlejoCorcho2024}, and a singular potential by Ianni, Le Coz and
Royer \cite{IanniLeCozRoyer2017}. Pelinovsky and Kevrekidis
\cite{PelinovskyKevrekidis2008} study persistence and spectral
stability under small decaying potentials, as well as particle
descriptions of dark-soliton motion. These results do not remove the
coordinate obstruction in the present hydrodynamic proof.

\subsection{The KdV structure and the two virial estimates}
Set
\begin{equation}\label{E:zeta}
\zeta_\beta=1-\hat\eta_\beta
=\sin^2\beta+\cos^2\beta\tanh^2(y\cos\beta),\qquad
\sigma=\begin{bmatrix}0&1\\1&0\end{bmatrix}.
\end{equation}
In hydrodynamic variables the Hamiltonian structure is
\[
\partial_t\begin{bmatrix}\eta\\\nu\end{bmatrix}
=JH_0'(\eta,\nu),\qquad J=-2\partial_x\sigma.
\]
The derivative in $J$ is the essential connection with KdV-type
equations. If $\mathcal L_\beta$ is the Hessian of the action at a
soliton, the linearized perturbation equation is $\partial_t=J\mathcal
L_\beta$. The transformed variable
\[
v=\sigma\mathcal L_\beta
\begin{bmatrix}\hat\theta\\\hat\xi\end{bmatrix}
\]
instead satisfies $v_t=-2\sigma\mathcal L_\beta v_y$ at frozen
parameter and linear order. Moving the Hessian to the other side of
the spatial derivative is the dual virial mechanism. It also
annihilates the linear translation forcing, because
$\mathcal L_\beta(\hat\eta_\beta',\hat\nu_\beta')^T=0$.

The local virial and rigidity method for asymptotic stability was
developed by Martel and Merle for generalized KdV
\cite{MartelMerle2000,MartelMerle2001,MartelMerle2005}.
The dual formulation and its improved virial coercivity are due to
Martel \cite{Martel2006}; see also Martel and Merle
\cite{MartelMerle2008}. For Gross--Pitaevskii, B\'ethuel, Gravejat and Smets
\cite[Section~1.3.3, equation~(34) and Proposition~7]{BGS_Asymp_Dark_2015}
applied the Hessian transform to the perturbation of a limiting
solution and used a virial with linear weight there. Their Propositions~8--9
and Corollary~1 combine coercive forms to give spacetime control on
that orbit. Their Proposition~2 also supplies localized momentum
monotonicity on the original orbit. A direct scheme
using separate virials for the evolving perturbation and a regularized
transform, followed by recovery of the original variables, was developed
by Kowalczyk, Martel and Mu\~noz for nonlinear Klein--Gordon
\cite{KMM2022} and by Martel for cubic--quintic Schr\"odinger
\cite[Section~3]{MartelCQ2022}. Here the weight in the transformed
Gross--Pitaevskii equation is $y\zeta_\beta$. There remains an important difference in dispersion: the
linearization at the constant background has the two branches
$\omega=\pm k\sqrt{k^2+4}$, so radiation travels in both directions.

For the homogeneous dual equation the weight $y\zeta_\beta$ gives
\[
\frac{d}{dt}\int y\zeta_\beta v_1v_2
=\langle\mathcal W_\beta v,v\rangle,\qquad
\mathcal W_\beta=\partial_y\mathcal L_\beta(y\zeta_\beta)
-(y\zeta_\beta)\mathcal L_\beta\partial_y.
\]
It is asymptotic to the usual linear weight at infinity and changes
the commutator in the soliton core. Theorem~\ref{T:virial-coercivity}
proves coercivity under the constraint transferred from modulation.
The proof is a scalar half-line estimate followed by a finite-dimensional
correction; the polynomial inequalities it uses are certified exactly
in Appendix~\ref{S:virial-certificate}. No numerical spectral hypothesis
is imposed. Related spectral questions for dual virials in the
Zakharov--Kuznetsov equation are discussed in
\cite{HolmerRoudenko2026}.

The unbounded weight and the two derivatives in $\mathcal L_\beta$
prevent this formal identity from being used directly at energy
regularity. We replace $y$ by $C\tanh(y/C)$ and regularize with
$(1-\gamma\partial_y^2)^{-1}$ after flattening the principal
coefficient. The sign is chosen so that the Fourier multiplier is
$(1+\gamma k^2)^{-1}$. Such regularized dual variables occur in
Kenig and Martel's work on Benjamin--Ono \cite{KenigMartel2009}.
Here a broad localization of the perturbation is made before the
dual transformation. This makes every nonlinear error proportional
to the same local energy that the primal virial controls.

The primal estimate controls local energy up to a localized density
term. The dual estimate controls the transformed variable, and an
elliptic inverse recovers that density term. A fixed positive linear
combination of the two bounded functionals is therefore increasing
at a rate comparable to the local energy. Integration gives
\eqref{E:main-estimate} directly, without extracting an asymptotic
orbit or invoking a nonlinear Liouville theorem.

\subsection{Why the black soliton is excluded}\label{S:black-intro}
The obstruction here is not instability of the black soliton:
its orbital and asymptotic stability are known, as described above.
It is the singularity of the variables and norms used in this proof.
At the center of a gray soliton,
\begin{equation}\label{E:black-core}
\zeta_\beta(0)=\sin^2\beta,\qquad
|q_\beta(0)|=|\sin\beta|,\qquad
\hat\nu_\beta(0)=\frac{\cos^2\beta}{\sin\beta}.
\end{equation}
The phase is undefined at the zero of $q_0$. Moreover the problem is
not confined to the value at one point. As $\beta\downarrow0$,
\[
\beta\hat\nu_\beta(\beta z)\longrightarrow\frac1{1+z^2},
\qquad
\|\hat\nu_\beta\|_2^2\sim\frac{\pi}{2\beta},
\qquad
\hat\nu_\beta\rightharpoonup\pi\delta_0
\quad\text{as measures}.
\]
These statements follow from the explicit formula: on the core scale
$\zeta_\beta(y)\sim\beta^2+y^2$, while
$\int\hat\nu_\beta=\pi-2\beta$ and the mass outside any neighborhood
of zero tends to zero. More explicitly, the substitution
$t=\tanh(y\cos\beta)$ gives
\[
\|\hat\nu_\beta\|_2^2
=\cos\beta+\frac{\cos(2\beta)}{\sin\beta}\left(\frac\pi2-\beta\right),
\qquad 0<\beta<\pi/2,
\]
which proves the stated $L^2$ asymptotic. Thus the gray profiles have
no limit in the hydrodynamic velocity space $L^2$.

Correspondingly, $\mathcal L_\beta$ contains inverse powers of
$\zeta_\beta$, the equivalence of physical and hydrodynamic energy
norms degenerates, and maintaining nonvanishing requires the density
perturbation to be small relative to $\sin^2\beta$.
The lower-right entry of $\mathcal W_\beta$ equals $\sin^4\beta$ at
zero. Testing its constrained form on $(0,f)^T$, with $f$ odd and
concentrated at the origin, shows that a coercivity constant in the
unweighted $\mathcal X$ norm cannot exceed $\sin^4\beta$.
A parameter-dependent conjugation can change this numerical scale,
but a conjugation that loses bounded invertibility cannot produce
uniform coercivity in the original norm.

In contrast, the complex profile $q_\beta$ and its physical energy
remain regular through $\beta=0$. A dynamics problem that drives
$\beta$ toward zero must therefore distinguish a failure of the
hydrodynamic chart from a physical loss of stability. The estimates
below can be used while $\beta$ remains in a compact gray interval.
Control up to, or across, zero would require coordinates that allow
vanishing amplitude, together with estimates uniform in the
corresponding physical norm. The frameworks of
\cite{GS_Asymp_Black_2015,KowalczykMartel2025} are relevant starting
points for that separate problem. The chosen real-valued hydrodynamic
momentum also records a phase-branch issue: its soliton values tend to
$\pi/2$ and $-\pi/2$ on the two sides of zero. Any continuation
through zero must account for this jump in the momentum coordinate.

\subsection{Organization}
Section~\ref{S:modulation} derives the modulation equations and
proves orbital control. Section~\ref{S:virial-static}
proves coercivity for the weight $y\zeta_\beta$.
Sections~\ref{S:regularization}--\ref{S:closure} construct the two
bounded virials and prove the spacetime and parameter estimates.
Section~\ref{S:asymptotics} first derives weak and time-averaged limits
from these estimates. The conserved energies of the unforced cubic
equation then yield the pointwise-in-time limit in the local energy norm.

\subsection{Disclosure on use of AI}

We formulated the approach taken in this paper, chose the dual virial
weight $y(1-\hat\eta_\beta(y))$, and worked out preliminary versions
of the Lyapunov and local virial estimates. Our motivation came from
the companion paper on effective dynamics \cite{HMEffectiveDynamics},
where tracking the soliton parameters requires quantitative control
of the radiation near the moving core. B\'ethuel, Gravejat and Smets
had already proved asymptotic stability of gray solitons using
compactness and a Liouville theorem for limiting solutions
\cite{BGS_Asymp_Dark_2015}. For our application, we sought to close
the spacetime estimates directly on the evolving perturbation over
each finite time interval.

We identified the constrained coercivity estimate needed for the dual
virial form and checked it numerically. ChatGPT Astra 6.2 found the
analytic proof of this estimate, combining a half-line argument with
an exact polynomial certificate. It also assisted us in carrying out
the spatial cutoff and frequency regularization arguments and in
drafting and revising the exposition. During revision, we used it to
work out how to combine local decay with the conserved energies of
Koch and Liao \cite{KochLiaoLow} to obtain strong convergence in the
local energy norm.

\section{Hydrodynamics, modulation, and orbital control}\label{S:modulation}
All functions in hydrodynamic coordinates are real. Inner products
are real $L^2$ pairings, extended by duality when derivatives are
distributions. Primes on profiles denote $y$ derivatives. Unless
stated otherwise, constants are uniform for $\beta\in I$.

\subsection{Hamiltonian formulation and the linearized operator}
The nonvanishing hydrodynamic phase space is the open set
\[
\{(\eta,\nu)\in\mathcal X:\sup\eta<1\}.
\]
On this set the energy and momentum are
\begin{equation}\label{E:conserved}
H_0(\eta,\nu)=\int\left\{
\frac{\eta_x^2}{8(1-\eta)}+\frac12(1-\eta)\nu^2+\frac12\eta^2
\right\}\,dx,\qquad P(\eta,\nu)=\frac12\int\eta\nu\,dx.
\end{equation}
Substitution of $u=\sqrt{1-\eta}\,e^{i\phi}$ in \eqref{E:GP}, with
$\phi_x=-\nu$, gives
\begin{equation}\label{E:hydro-system}
\begin{aligned}
\eta_t&=-2\partial_x((1-\eta)\nu),\\
\nu_t&=-2\partial_x\left(
-\frac{\eta_{xx}}{4(1-\eta)}
-\frac{\eta_x^2}{8(1-\eta)^2}-\frac12\nu^2+\eta\right).
\end{aligned}
\end{equation}
These are precisely $\partial_t(\eta,\nu)^T=JH_0'(\eta,\nu)$, since
\[
H_0'(\eta,\nu)=\begin{bmatrix}
-\frac14\partial_x((1-\eta)^{-1}\eta_x)
+\frac18(1-\eta)^{-2}\eta_x^2-\frac12\nu^2+\eta\\
(1-\eta)\nu
\end{bmatrix},\qquad
P'(\eta,\nu)=\frac12\begin{bmatrix}\nu\\\eta\end{bmatrix}.
\]
We use Mohamad's hydrodynamic well-posedness theorem
\cite[Theorem~1.2]{Mohamad2014}, also stated in
\cite[Theorem~1]{BGSChains2014}. It gives continuous dependence in
$H^1\times L^2$, conservation of \eqref{E:conserved}, and continuation
as long as $1-\eta$ stays bounded below by a positive constant.
The corresponding higher-regularity theory applies to smooth data.
For the Madelung map on the nonvanishing quotient space see
\cite[Theorem~1.6]{Wegner2023}.
The normalization of \cite{BGS_Asymp_Dark_2015} is converted to ours by
$u(t,x)=\Psi(2t,\sqrt2x)$. Reconstruction determines the solution up
to a constant phase.

Define the action $Z_\beta=H_0-2\sin\beta\,P$. The soliton satisfies
\begin{equation}\label{E:Zp-soliton}
Z_\beta'(\hat\eta_\beta,\hat\nu_\beta)=0.
\end{equation}
Its Hessian is
\begin{equation}\label{E:L-definition}
\begin{aligned}
\mathcal L_\beta&=
\begin{bmatrix}\mathcal A_\beta&-\sin\beta/\zeta_\beta\\
-\sin\beta/\zeta_\beta&\zeta_\beta\end{bmatrix},\\
\mathcal A_\beta&=-\frac14\partial_y\zeta_\beta^{-1}\partial_y
+1+\frac{\zeta_\beta''}{4\zeta_\beta^2}
-\frac{(\zeta_\beta')^2}{4\zeta_\beta^3}.
\end{aligned}
\end{equation}
In particular the off-diagonal entry is
$-\hat\nu_\beta-\sin\beta=-\sin\beta/\zeta_\beta$.
Differentiating \eqref{E:Zp-soliton} in translation and in $\beta$ gives
\begin{equation}\label{E:L-identities}
\mathcal L_\beta\begin{bmatrix}\hat\eta_\beta'\\\hat\nu_\beta'\end{bmatrix}=0,
\qquad
\mathcal L_\beta\partial_\beta
\begin{bmatrix}\hat\eta_\beta\\\hat\nu_\beta\end{bmatrix}
=\cos\beta\,\sigma
\begin{bmatrix}\hat\eta_\beta\\\hat\nu_\beta\end{bmatrix}.
\end{equation}
The profiles are even, and all their derivatives in $y$ and $\beta$
decay exponentially, uniformly on $I$.

\begin{lemma}[Hessian spectrum and constrained energy]\label{L:L-spectrum}
The self-adjoint operator $\mathcal L_\beta$ on $L^2\times L^2$, with
domain $H^2\times L^2$, has one simple negative eigenvalue,
kernel $\spn\{(\hat\eta_\beta',\hat\nu_\beta')^T\}$, and a positive
spectral gap after these directions are removed. If $F\in\mathcal X$
satisfies
\[
\left\langle F,\sigma\begin{bmatrix}\hat\eta_\beta\\\hat\nu_\beta\end{bmatrix}
\right\rangle=0,\qquad
\left\langle F,\begin{bmatrix}\hat\eta_\beta'\\\hat\nu_\beta'\end{bmatrix}
\right\rangle=0,
\]
then, for a constant $c_I>0$,
\begin{equation}\label{E:energy-coercivity}
\langle\mathcal L_\beta F,F\rangle\geq c_I\|F\|_{\mathcal X}^2.
\end{equation}
\end{lemma}
\begin{proof}
Completing the square gives
\begin{equation}\label{E:L-square}
\langle\mathcal L_\beta F,F\rangle
=\langle H_\beta F_1,F_1\rangle
+\int\zeta_\beta(F_2-\sin\beta\,\zeta_\beta^{-2}F_1)^2,
\qquad H_\beta=\mathcal A_\beta-\sin^2\beta\,\zeta_\beta^{-3}.
\end{equation}
The scalar operator $H_\beta$ is uniformly elliptic and approaches
$-\frac14\partial_y^2+\cos^2\beta$ at infinity. By
\eqref{E:L-identities} and
$\hat\nu_\beta'=\sin\beta\,\zeta_\beta^{-2}\hat\eta_\beta'$,
its kernel contains $\hat\eta_\beta'$, which has exactly one zero.
Sturm oscillation gives one negative eigenvalue and a simple zero,
with the remaining spectrum bounded away from zero.
The invertible triangular change in \eqref{E:L-square} preserves the
negative index and identifies the kernel of $\mathcal L_\beta$.
Since $\zeta_\beta\geq\sin^2\beta$, it also gives a positive
essential gap and the stated self-adjoint realization.

For clarity, we verify the sign that removes the negative direction.
Direct integration of \eqref{E:hydro-profile} gives, for $\beta>0$,
\begin{equation}\label{E:profile-momentum}
P(\hat\eta_\beta,\hat\nu_\beta)
=\frac\pi2-\beta-\sin\beta\cos\beta,
\qquad \partial_\beta P=-2\cos^2\beta.
\end{equation}
For example, $\hat\eta_\beta\hat\nu_\beta
=\hat\nu_\beta-\sin\beta\hat\eta_\beta$,
$\int\hat\nu_\beta=\pi-2\beta$, and
$\int\hat\eta_\beta=2\cos\beta$.
The even vector $\sec\beta\,\partial_\beta
(\hat\eta_\beta,\hat\nu_\beta)^T$ therefore has negative form value
$-4\cos\beta$ and is form-orthogonal to every $F$ obeying the first
constraint. Since the negative index is one, the form is nonnegative
on that constraint space. A null vector of the restricted form
satisfies
$\mathcal L_\beta F=\lambda\sigma(\hat\eta_\beta,\hat\nu_\beta)^T$.
Pairing with $\sec\beta\,\partial_\beta
(\hat\eta_\beta,\hat\nu_\beta)^T$ gives $-4\lambda\cos\beta=0$.
Thus its only null direction is the translation kernel.
The positive essential gap now gives a strict $L^2\times L^2$ bound
under both constraints. One may see this directly by decomposing
into the negative eigenvector, the kernel, and the positive spectral
space: the negative form value just computed makes the
Cauchy--Schwarz inequality for the constrained negative coefficient
strict. Finally, the positive derivative coefficient in
\eqref{E:L-definition} upgrades the bound to \eqref{E:energy-coercivity}.
Continuity of the coefficients and a finite cover of $I$ give uniformity.
\end{proof}

\subsection{Modulation and the exact nonlinear equation}
We impose
\begin{equation}\label{E:orth}
\left\langle\begin{bmatrix}\hat\theta\\\hat\xi\end{bmatrix},
\sigma\begin{bmatrix}\hat\eta_\beta\\\hat\nu_\beta\end{bmatrix}\right\rangle=0,
\qquad
\left\langle\begin{bmatrix}\hat\theta\\\hat\xi\end{bmatrix},
\begin{bmatrix}\hat\eta_\beta'\\\hat\nu_\beta'\end{bmatrix}\right\rangle=0.
\end{equation}
We use the momentum and translation constraints in the form of
\cite[Theorem~3, equation~(6)]{BGS_Asymp_Dark_2015}. The first condition
controls the negative energy direction, and the second fixes translation.
In their normalization, estimate~(9) gives a quadratic bound on the speed
derivative in the full energy norm. In Lemma~\ref{L:modulation-odes}
we keep the spatial weight and identify the linear term in the center
equation.

\begin{lemma}[Modulation]\label{L:modulation}
In a sufficiently small energy neighborhood of a translated gray
soliton, \eqref{E:orth} determines a unique nearby pair $(\beta,a)$,
depending $C^1$ on the data. Along an energy solution these parameters
are $C^1$ in time.
\end{lemma}
\begin{proof}
At a soliton the parameter Jacobian is diagonal by parity. Its
diagonal entries, up to the translation sign convention, are
\[
\left\langle\partial_\beta
\begin{bmatrix}\hat\eta_\beta\\\hat\nu_\beta\end{bmatrix},
\sigma\begin{bmatrix}\hat\eta_\beta\\\hat\nu_\beta\end{bmatrix}\right\rangle
=-4\cos^2\beta,
\qquad
\left\|\begin{bmatrix}\hat\eta_\beta'\\\hat\nu_\beta'\end{bmatrix}\right\|_2^2>0.
\]
The implicit function theorem applies uniformly on $I$. The pairings
are differentiable with respect to translations even at energy
regularity, since derivatives can be transferred to the smooth tests.
For the time dependence, write each constraint in laboratory
coordinates as a pairing of the solution with a translated, smooth
profile test, minus its soliton value. Holding $a,\beta$ fixed, its
time derivative is obtained from \eqref{E:hydro-system} and is
continuous since the right side belongs continuously to
$H^{-1}\times H^{-2}$. The implicit function theorem in
$(t,a,\beta)$ thus gives $C^1$ time dependence.
\end{proof}

The perturbation obeys the exact identity
\begin{equation}\label{E:pert-eqn}
\begin{aligned}
\partial_t\begin{bmatrix}\hat\theta\\\hat\xi\end{bmatrix}
={}&J\mathcal L_\beta\begin{bmatrix}\hat\theta\\\hat\xi\end{bmatrix}
-\dot\beta\,\partial_\beta
\begin{bmatrix}\hat\eta_\beta\\\hat\nu_\beta\end{bmatrix}
+(\dot a-2\sin\beta)
\begin{bmatrix}\hat\eta_\beta'\\\hat\nu_\beta'\end{bmatrix}\\
&+(\dot a-2\sin\beta)\partial_y
\begin{bmatrix}\hat\theta\\\hat\xi\end{bmatrix}
+J\mathcal N_\beta.
\end{aligned}
\end{equation}
Here $\mathcal N_\beta$ is the nonlinear Taylor remainder of $H_0'$ at
$(\hat\eta_\beta,\hat\nu_\beta)^T$. There is no nonlinear remainder
from momentum. The identities $JP'=-(\eta,\nu)_y^T$ and
$JP''=-\partial_y$ account for the speed $2\sin\beta$ in
\eqref{E:pert-eqn}.

We record two complementary descriptions of $\mathcal N_\beta$.
First, it is the variational derivative in $(\hat\theta,\hat\xi)$ of
the integral of the exact cubic energy remainder
\begin{equation}\label{E:cubic-density}
\mathcal R_\beta(y,\hat\theta,\hat\theta_y,\hat\xi)
=\frac{\hat\theta}{8\zeta_\beta(\zeta_\beta-\hat\theta)}
\left(\hat\theta_y+\frac{\hat\eta_\beta'}{\zeta_\beta}\hat\theta\right)^2
-\frac12\hat\theta\hat\xi^2.
\end{equation}
This follows by subtracting the Taylor polynomial of degree two from
the integrand in \eqref{E:conserved}; the kinetic remainder is a square.
It implies
$|\int\mathcal R_\beta|\leq C\|\hat\theta\|_{H^1}
\| (\hat\theta,\hat\xi)^T\|_{\mathcal X}^2$.
Second, the distributional form needed for regularization is
\begin{equation}\label{E:remainder-GF}
(\mathcal N_\beta)_1=\partial_yG+F,
\qquad (\mathcal N_\beta)_2=-\hat\theta\hat\xi,
\end{equation}
where
\begin{align*}
G={}&-\frac{\hat\theta}{4\zeta_\beta(\zeta_\beta-\hat\theta)}
\left(\hat\theta_y+\frac{\hat\eta_\beta'}{\zeta_\beta}\hat\theta\right),\\
F={}&\frac18\left[
\left(\frac1{(\zeta_\beta-\hat\theta)^2}
-\frac1{\zeta_\beta^2}-\frac{2\hat\theta}{\zeta_\beta^3}\right)
(\hat\eta_\beta')^2\right.\\
&\left.\hspace{12mm}
+2\left(\frac1{(\zeta_\beta-\hat\theta)^2}
-\frac1{\zeta_\beta^2}\right)\hat\eta_\beta'\hat\theta_y
+\frac{\hat\theta_y^2}{(\zeta_\beta-\hat\theta)^2}\right]
-\frac12\hat\xi^2.
\end{align*}
On the small nonvanishing energy neighborhood, $G\in L^2$ and
$F\in L^1\subset H^{-1}$. All nonlinear pairings with profile tests
are therefore legitimate. Both descriptions will be used below:
\eqref{E:cubic-density} exposes the cancellation in the primal virial,
whereas \eqref{E:remainder-GF} gives the low-regularity source bounds.

\begin{lemma}[Parameter equations]\label{L:modulation-odes}
Choose $k>0$ sufficiently small in terms of $I$ and let
$\rho(y)=\sech(ky)$. On an orbitally small solution,
\begin{equation}\label{E:modulation-bounds}
\begin{aligned}
|\dot\beta|&\leq C\int\rho^2
(\hat\theta_y^2+\hat\theta^2+\hat\xi^2),\\
\left|\dot a-2\sin\beta+
\frac{2\left\langle\begin{bmatrix}\hat\theta\\\hat\xi\end{bmatrix},
\mathcal L_\beta\sigma
\begin{bmatrix}\hat\eta_\beta''\\\hat\nu_\beta''\end{bmatrix}\right\rangle}
{\left\|\begin{bmatrix}\hat\eta_\beta'\\\hat\nu_\beta'\end{bmatrix}\right\|_2^2}
\right|&\leq C\int\rho^2
(\hat\theta_y^2+\hat\theta^2+\hat\xi^2).
\end{aligned}
\end{equation}
In particular $|\dot a-2\sin\beta|\leq C\delta$ if the perturbation
has energy norm at most $C\delta$.
\end{lemma}
\begin{proof}
Differentiate \eqref{E:orth} and use \eqref{E:pert-eqn}. The linear
term in the first equation vanishes, because
\[
\left\langle J\mathcal L_\beta
\begin{bmatrix}\hat\theta\\\hat\xi\end{bmatrix},
\sigma\begin{bmatrix}\hat\eta_\beta\\\hat\nu_\beta\end{bmatrix}\right\rangle
=2\left\langle\begin{bmatrix}\hat\theta\\\hat\xi\end{bmatrix},
\mathcal L_\beta
\begin{bmatrix}\hat\eta_\beta'\\\hat\nu_\beta'\end{bmatrix}\right\rangle=0.
\]
The pure profile cross terms vanish by parity. Integrating derivatives
onto the tests gives the exact two-equation system
\begin{equation}\label{E:modulation-system}
\left\{\begin{aligned}
&\left(-4\cos^2\beta-
\left\langle\begin{bmatrix}\hat\theta\\\hat\xi\end{bmatrix},
\sigma\partial_\beta\begin{bmatrix}\hat\eta_\beta\\\hat\nu_\beta\end{bmatrix}
\right\rangle\right)\dot\beta\\
&\quad+\left\langle\begin{bmatrix}\hat\theta\\\hat\xi\end{bmatrix},
\sigma\begin{bmatrix}\hat\eta_\beta'\\\hat\nu_\beta'\end{bmatrix}\right\rangle
(\dot a-2\sin\beta)
=2\left\langle\mathcal N_\beta,
\begin{bmatrix}\hat\eta_\beta'\\\hat\nu_\beta'\end{bmatrix}\right\rangle,\\[2pt]
&\left\langle\begin{bmatrix}\hat\theta\\\hat\xi\end{bmatrix},
\partial_\beta\begin{bmatrix}\hat\eta_\beta'\\\hat\nu_\beta'\end{bmatrix}
\right\rangle\dot\beta\\
&\quad+\left(\left\|\begin{bmatrix}\hat\eta_\beta'\\\hat\nu_\beta'\end{bmatrix}\right\|_2^2
-\left\langle\begin{bmatrix}\hat\theta\\\hat\xi\end{bmatrix},
\begin{bmatrix}\hat\eta_\beta''\\\hat\nu_\beta''\end{bmatrix}\right\rangle\right)
(\dot a-2\sin\beta)\\
&\qquad=-2\left\langle\mathcal L_\beta
\begin{bmatrix}\hat\theta\\\hat\xi\end{bmatrix}+\mathcal N_\beta,
\sigma\begin{bmatrix}\hat\eta_\beta''\\\hat\nu_\beta''\end{bmatrix}\right\rangle.
\end{aligned}\right.
\end{equation}
The matrix is a small perturbation of an invertible diagonal matrix.
After integrating $G'$ once, \eqref{E:remainder-GF} bounds each
nonlinear pairing by $C\int\rho^2
(\hat\theta_y^2+\hat\theta^2+\hat\xi^2)$.
The linear pairing in the second equation has size at most
$C\|\rho(\hat\theta,\hat\xi)^T\|_2$, since $\mathcal L_\beta$ can
be moved to the smooth test. The off-diagonal coefficients have the
same localized linear bound. Inverting the matrix gives first the
linear bound for $\dot a-2\sin\beta$ and then the quadratic bound
for $\dot\beta$; substituting back gives the second estimate in
\eqref{E:modulation-bounds}. The pairings in this system are continuous at energy regularity.
\end{proof}

\subsection{Global orbital control}
\begin{proposition}\label{P:orbital}
For data in Theorem~\ref{T:main}, the modulation extends to all real
times, and
\begin{equation}\label{E:orbital-control}
\sup_t\left\|\begin{bmatrix}\hat\theta(t)\\\hat\xi(t)\end{bmatrix}\right\|_{\mathcal X}
\leq C\delta,\qquad
\sup_t|\beta(t)-\beta(0)|\leq C\delta^2,
\qquad
1-\eta\geq\frac12\min_{\beta\in I}\sin^2\beta.
\end{equation}
\end{proposition}
\begin{proof}
Write $\beta_* =\beta(0)$ for the fitted initial parameter. It differs
from the reference $\beta_0$ by $O(\delta)$.
The first orthogonality condition and the quadratic nature of momentum
give the exact identity
\[
P(\eta,\nu)=P(\hat\eta_\beta,\hat\nu_\beta)
+\frac12\int\hat\theta\hat\xi.
\]
Since $\partial_\beta P(\hat\eta_\beta,\hat\nu_\beta)
=-2\cos^2\beta$ is bounded away from zero on $I$, conservation yields
\begin{equation}\label{E:beta-est2}
|\beta(t)-\beta_*|
\leq C\left(\delta^2+
\left\|\begin{bmatrix}\hat\theta(t)\\\hat\xi(t)\end{bmatrix}\right\|_{\mathcal X}^2\right).
\end{equation}
To compare energy at different parameters, subtract the soliton value
from the action. Its parameter derivative at any fixed datum is
\[
\partial_\beta\{Z_\beta(\eta,\nu)
-Z_\beta(\hat\eta_\beta,\hat\nu_\beta)\}
=-2\cos\beta\{P(\eta,\nu)-P(\hat\eta_\beta,\hat\nu_\beta)\}.
\]
At the initial datum and $\beta=\beta_*$ this is $O(\delta^2)$,
and its next parameter derivative is bounded. Conservation of both
$H_0$ and $P$ and Taylor's theorem in $\beta$ therefore give
\begin{align*}
&Z_{\beta(t)}(\eta(t),\nu(t))
-Z_{\beta(t)}(\hat\eta_{\beta(t)},\hat\nu_{\beta(t)})\\
&\qquad\leq C\delta^2+C\delta^2|\beta(t)-\beta_*|
+C|\beta(t)-\beta_*|^2.
\end{align*}
On the other hand, the exact expansion \eqref{E:cubic-density} and
\eqref{E:energy-coercivity} bound the left side below by a positive
multiple of $\| (\hat\theta,\hat\xi)^T\|_{\mathcal X}^2$ on a small
energy neighborhood. Hence
\[
\left\|\begin{bmatrix}\hat\theta(t)\\\hat\xi(t)\end{bmatrix}\right\|_{\mathcal X}^2
\leq C\delta^2+C
\left\|\begin{bmatrix}\hat\theta(t)\\\hat\xi(t)\end{bmatrix}\right\|_{\mathcal X}^4.
\]
A continuity argument with a fixed small bootstrap radius gives the
first bound in \eqref{E:orbital-control}, and \eqref{E:beta-est2}
gives the second. Sobolev embedding and
$\zeta_\beta\geq\min_I\sin^2\beta$ give the last bound.
These estimates keep the solution strictly inside the modulation
neighborhood and $\beta$ strictly inside $I$. Energy-space
well-posedness and Lemma~\ref{L:modulation} continue the decomposition
at every finite endpoint. The same argument backwards in time
completes the proof. No spacetime estimate is used at this stage.
\end{proof}

\section{The dual virial and its coercivity}\label{S:virial-static}
The orbital estimate controls the size of the perturbation but does
not force radiation to leave the soliton. We now identify a quadratic
form that supplies this missing spacetime control.

\subsection{The dual equation and the weight}
Temporarily fix $\beta$ and take a smooth, rapidly decreasing solution
of the homogeneous linearized equation. Define
\begin{equation}\label{E:formal-dual}
v=\sigma\mathcal L_\beta
\begin{bmatrix}\hat\theta\\\hat\xi\end{bmatrix}.
\end{equation}
Because $J=-2\sigma\partial_y$, its evolution is
\[
v_t=-2\sigma\mathcal L_\beta\partial_y v.
\]
This identifies the principal quadratic form. Section~\ref{S:dual}
retains the nonlinear terms, parameter derivatives, and regularization
errors in the full equation.
For any real weight $m$ for which the pairings are meaningful,
\begin{equation}\label{E:dual-commutator}
\frac{d}{dt}\int m v_1v_2
=-2\langle m\mathcal L_\beta v_y,v\rangle
=\langle(\partial_y\mathcal L_\beta m-m\mathcal L_\beta\partial_y)v,v\rangle.
\end{equation}
The products in the last expression are compositions of operators;
in particular $\partial_y$ differentiates everything to its right.

For completeness, the scalar identity behind the commutator is
\begin{equation}\label{E:scalar-commutator}
\partial_y(-\partial_y p\partial_y+b)m
-m(-\partial_y p\partial_y+b)\partial_y
=-\partial_y(3pm'+p'm)\partial_y+(bm)'-(pm')''.
\end{equation}
It follows by using $\partial_y m=m\partial_y+m'$ twice and
collecting the zeroth and second order terms. The first order
coefficient is the derivative of the second order coefficient, as
required by self-adjointness.
Apply this identity with $p=(4\zeta_\beta)^{-1}$ and with the
zeroth order coefficient of $\mathcal A_\beta$ in place of $b$.
For $m=y\zeta_\beta$, the off-diagonal entry becomes the constant
$-\sin\beta$, and \eqref{E:dual-commutator} yields the operator

\begin{equation}\label{E:W-compact}
\mathcal W_\beta=
\begin{bmatrix}
-\partial_y(\frac34+\frac{y\zeta_\beta '}{2\zeta_\beta })\partial_y
+\zeta_\beta +y\zeta_\beta '-\frac14(\log\zeta_\beta )''&-\sin\beta \\
-\sin\beta &\zeta_\beta ^2+2y\zeta_\beta \zeta_\beta '
\end{bmatrix}.
\end{equation}
Its coefficients are smooth on a gray branch, and its quadratic form
is well-defined on $H^1\times L^2$. Define
\[
h_\beta=\sigma\partial_\beta
\begin{bmatrix}\hat\eta_\beta\\\hat\nu_\beta\end{bmatrix}
=\begin{bmatrix}\partial_\beta\hat\nu_\beta\\\partial_\beta\hat\eta_\beta\end{bmatrix}.
\]

\begin{theorem}[Virial coercivity]\label{T:virial-coercivity}
For each $0<|\beta|<\pi/2$ there is $\kappa_\beta>0$ such that
\begin{equation}\label{E:virial-coercivity}
\left\langle\mathcal W_\beta F,F\right\rangle
\geq\kappa_\beta\bigl(\left\|F_1\right\|_{H^1}^2+\left\|F_2\right\|_2^2\bigr)
\end{equation}
for every real $F\in\mathcal X$ with $\langle F,h_\beta\rangle=0$.
The constant can be chosen uniformly on each compact subinterval
of either gray branch; we denote this constant by $\kappa_I$ on $I$.
\end{theorem}

We prove the estimate first for $0<\beta<\pi/2$. In the proof, write $\mathfrak w_\beta(F,G)$ for the quadratic-form pairing of $\mathcal W_\beta$
and $\mathfrak w_\beta(F)=\mathfrak w_\beta(F,F)$. Denote its coefficient functions by
\begin{equation}\label{E:virial-coefficients}
\begin{aligned}
a_\beta&=\frac34+\frac{y\zeta_\beta'}{2\zeta_\beta},&
b_\beta&=\zeta_\beta+y\zeta_\beta'-\frac14(\log\zeta_\beta)'',\\
d_\beta&=\zeta_\beta^2+2y\zeta_\beta\zeta_\beta'.&&
\end{aligned}
\end{equation}
Since $y\zeta_\beta'\geq0$, we have $a_\beta\geq3/4$ and
$d_\beta\geq\sin^4\beta$. Completing the square gives
\begin{equation}\label{E:virial-square}
\begin{aligned}
K_\beta&=-\partial_y a_\beta\partial_y+b_\beta
-\frac{\sin^2\beta}{d_\beta},\\
\mathfrak w_\beta(F)&=\langle K_\beta F_1,F_1\rangle
+\int d_\beta\left(F_2-\frac{\sin\beta}{d_\beta}F_1\right)^2\,dy.
\end{aligned}
\end{equation}
To pass from a trace condition at the origin to the modulation
constraint, compare $h_\beta$ with the image of the scaled profile:
\begin{equation}\label{E:virial-comparison-vector}
r_\beta=h_\beta-\frac{2}{\cos^3\beta}\mathcal W_\beta
\begin{bmatrix}\hat\eta_\beta\\\hat\nu_\beta\end{bmatrix}.
\end{equation}
The two identities needed for this comparison are
\begin{equation}\label{E:virial-profile-identities}
\mathfrak w_\beta\!\left(
\begin{bmatrix}\hat\eta_\beta\\\hat\nu_\beta\end{bmatrix}\right)=0,
\qquad
\left\langle h_\beta,
\begin{bmatrix}\hat\eta_\beta\\\hat\nu_\beta\end{bmatrix}\right\rangle
=-4\cos^2\beta.
\end{equation}
For the first, use the commutator identity and
$\mathcal L_\beta(\hat\eta_\beta',\hat\nu_\beta')^T=0$:
\[
\mathfrak w_\beta\!\left(
\begin{bmatrix}\hat\eta_\beta\\\hat\nu_\beta\end{bmatrix}\right)
=-2\left\langle\mathcal L_\beta
\begin{bmatrix}\hat\eta_\beta'\\\hat\nu_\beta'\end{bmatrix},
y\zeta_\beta\begin{bmatrix}\hat\eta_\beta\\\hat\nu_\beta\end{bmatrix}
\right\rangle=0.
\]
The second follows by differentiating
$\int\hat\eta_\beta\hat\nu_\beta
=\pi-2\beta-2\sin\beta\cos\beta$.

\begin{lemma}[Half-line inequalities]\label{L:virial-halfline}
For $y>0$ define
\begin{equation}\label{E:virial-test-function}
\psi_\beta(y)=
\frac{\tanh(y\cos\beta)\,\sech^{7/8}(y\cos\beta)}
{\bigl(\sin^2\beta+3\cos^2\beta\tanh^2(y\cos\beta)\bigr)^{1/3}},
\qquad
\Lambda_\beta(y)=\frac{K_\beta\psi_\beta(y)}{\psi_\beta(y)}.
\end{equation}
Then
\begin{equation}\label{E:virial-pointwise}
\begin{aligned}
\Lambda_\beta(y)&\geq\tfrac14\cos^2\beta,\\
\frac{\bigl((r_\beta)_1+\sin\beta\,(r_\beta)_2/d_\beta\bigr)^2}
{\Lambda_\beta}
+\frac{(r_\beta)_2^2}{d_\beta}
&\leq\frac{15}{2}\sech^2(y\cos\beta).
\end{aligned}
\end{equation}
\end{lemma}
\begin{proof}
Appendix~\ref{S:virial-certificate} reduces these inequalities to
nonnegativity of two explicitly given polynomials and verifies their
Bernstein coefficients by exact integer arithmetic.
\end{proof}

Extend $\Lambda_\beta$ evenly to the line. If $G_1(0)=0$, the
positive-function identity on the two half-lines and
\eqref{E:virial-square} give
\begin{equation}\label{E:virial-trace-bound}
\mathfrak w_\beta(G)\geq\int\left[
\Lambda_\beta G_1^2+d_\beta
\left(G_2-\frac{\sin\beta}{d_\beta}G_1\right)^2\right]\,dy.
\end{equation}
The omitted term on $y>0$ is
$\int_0^\infty a_\beta\psi_\beta^2|(G_1/\psi_\beta)'|^2\,dy$;
on the other half-line reflect $\psi_\beta$.
First apply the identity to compactly supported functions away from zero,
then use density in the trace-zero space.
For fixed $\beta$, $\Lambda_\beta$ is bounded at both endpoints:
$\psi_\beta$ has a simple odd zero at the origin and decays
exponentially at infinity. By Cauchy--Schwarz and
\eqref{E:virial-pointwise},
\begin{equation}\label{E:virial-residual-bound}
|\langle r_\beta,G\rangle|^2
\leq\frac{15}{\cos\beta}\,\mathfrak w_\beta(G)
\qquad\text{if }G_1(0)=0.
\end{equation}
Here one uses
\[
r_\beta\cdot G=
\left((r_\beta)_1+\frac{\sin\beta}{d_\beta}(r_\beta)_2\right)G_1
+(r_\beta)_2\left(G_2-\frac{\sin\beta}{d_\beta}G_1\right)
\]
and $\int\sech^2(y\cos\beta)\,dy=2/\cos\beta$.
Moreover, \eqref{E:virial-trace-bound} controls
$\|G\|_{L^2\times L^2}^2$ for fixed $\beta$.
Returning to \eqref{E:W-compact}, $a_\beta\geq3/4$ and the bounded
zeroth-order coefficients then give
\begin{equation}\label{E:virial-trace-coercivity}
\mathfrak w_\beta(G)\gtrsim_\beta\|G\|_{\mathcal X}^2
\qquad\text{if }G_1(0)=0.
\end{equation}

\begin{proof}[Proof of Theorem~\ref{T:virial-coercivity}]
Given $F\perp h_\beta$, set $k=\frac12\cos\beta\,F_1(0)$ and write
\[
F=G+\frac{2k}{\cos^3\beta}
\begin{bmatrix}\hat\eta_\beta\\\hat\nu_\beta\end{bmatrix}.
\]
Since $\hat\eta_\beta(0)=\cos^2\beta$, this makes $G_1(0)=0$.
By \eqref{E:virial-profile-identities},
$\langle h_\beta,G\rangle=8k/\cos\beta$. Therefore
\begin{align}
\mathfrak w_\beta(F)&=\mathfrak w_\beta(G)+\frac{16}{\cos\beta}k^2
-2k\langle r_\beta,G\rangle\nonumber\\
&\geq\left(1-\sqrt{15/16}\right)
\left(\mathfrak w_\beta(G)+\frac{16}{\cos\beta}k^2\right).
\label{E:virial-strict-margin}
\end{align}
The last step uses \eqref{E:virial-residual-bound} and $2AB\leq A^2+B^2$.
Together with \eqref{E:virial-trace-coercivity} and the finite $\mathcal X$ norm of the profile,
this proves \eqref{E:virial-coercivity}.
All quantities used to pass to the $\mathcal X$ norm are bounded
uniformly on compact $I\Subset(0,\pi/2)$, by the explicit profiles,
$d_\beta\geq\sin^4\beta$, and \eqref{E:virial-pointwise}.
Changing the sign of the second component intertwines the operators
at $\beta$ and $-\beta$ and takes $h_\beta$ to $h_{-\beta}$,
which proves the statement on the other gray branch.
\end{proof}

\begin{corollary}[The virial spectrum]\label{C:virial-spectrum}
For every fixed gray parameter, $\mathcal W_\beta$ has one simple
negative eigenvalue, trivial kernel, and a positive gap above zero.
Moreover,
\[
\langle\mathcal W_\beta^{-1}h_\beta,h_\beta\rangle<0.
\]
\end{corollary}
\begin{proof}
The form is closed and semibounded on $\mathcal X$.
Theorem~\ref{T:virial-coercivity} bounds it below by
$\kappa_\beta\|F\|_2^2$ on a codimension-one subspace.
Thus the spectral subspace below $\kappa_\beta$ has dimension at most one.
The profile has zero form value but is not a kernel vector:
the second component of
$\mathcal W_\beta(\hat\eta_\beta,\hat\nu_\beta)^T$ at zero is
$-\sin\beta\cos^4\beta\ne0$.
Perturbing it in a direction with a nonzero form pairing gives a negative
vector, proving the spectral assertions.
If $u=\mathcal W_\beta^{-1}h_\beta$, then $\mathfrak w_\beta(u,G)=0$
for $G\perp h_\beta$. The value $\mathfrak w_\beta(u)=\langle u,h_\beta\rangle$
cannot vanish by the coercivity theorem; if it were positive,
the form would be nonnegative on
$h_\beta^\perp\oplus\operatorname{span}\{u\}=\mathcal X$,
which is also impossible.
\end{proof}

\begin{corollary}[The transferred constraint]\label{C:transferred-form}
If the perturbation satisfies the first condition in \eqref{E:orth},
then its formal dual variable \eqref{E:formal-dual} satisfies
$\langle v,h_\beta\rangle=0$. Whenever $v\in\mathcal X$,
\[
\langle\mathcal W_\beta v,v\rangle\geq\kappa_I\|v\|_{\mathcal X}^2.
\]
\end{corollary}
\begin{proof}
Self-adjointness and \eqref{E:L-identities} give
\[
\left\langle\sigma\mathcal L_\beta
\begin{bmatrix}\hat\theta\\\hat\xi\end{bmatrix},
\sigma\partial_\beta\begin{bmatrix}\hat\eta_\beta\\\hat\nu_\beta\end{bmatrix}
\right\rangle
=\cos\beta\left\langle\begin{bmatrix}\hat\theta\\\hat\xi\end{bmatrix},
\sigma\begin{bmatrix}\hat\eta_\beta\\\hat\nu_\beta\end{bmatrix}\right\rangle=0.
\]
Apply Theorem~\ref{T:virial-coercivity}. The next section replaces the
formal dual variable by one defined on every energy perturbation,
and quantifies the resulting defect in this constraint.
\end{proof}

\section{Truncation, regularization, and recovery of the density}
\label{S:regularization}
We now work with the global small solution of Proposition~\ref{P:orbital}.
The purpose of this section is to define the dual variable at energy
regularity and connect it back to the original perturbation. The
parameters $A,B,\gamma$ below are fixed in time. We choose them in the order
\begin{equation}\label{E:choice-order}
B\text{ large},\qquad \gamma\text{ small},\qquad
A\text{ large depending on }B,\gamma,\qquad \delta\text{ small}.
\end{equation}
The broad radius $A$ localizes the perturbation before the dual
transformation; the radius $B$ truncates the dual virial weight.

\subsection{Weights and the smoothing operator}
For $C\geq1$ define
\begin{equation}\label{E:weights}
\chi_C(y)=\sech(y/C),\qquad
\phi_C(y)=C\tanh(y/C),\qquad \phi_C'=\chi_C^2.
\end{equation}
Thus the truncated weight is $\zeta_\beta\phi_C$. Its size is bounded
by $C$, and it tends to $y\zeta_\beta$ on bounded intervals as
$C\to\infty$. For a test vector $F$ put
\[
D_C(F)=\int\chi_C^2(|F_1'|^2+|F_1|^2+|F_2|^2).
\]
This is uniformly equivalent to $\|\chi_CF\|_{\mathcal X}^2$ for
$C\geq1$. For the evolving perturbation write
\begin{equation}\label{E:local-energy}
\mathcal D_A(t)=\int\chi_A^2
(\hat\theta_y^2+\hat\theta^2+\hat\xi^2).
\end{equation}
The fixed core weight $\rho=\sech(ky)$ from
Lemma~\ref{L:modulation-odes} is chosen to decay more slowly than
the profile coefficients. In particular $\rho\lesssim\chi_A$ for
all sufficiently large $A$.

\begin{lemma}[Resolvent estimates]\label{L:resolvent}
For $0<\gamma\leq1$, let
\[
R_\gamma=(1-\gamma\partial_y^2)^{-1},\qquad
R_\gamma f=\frac1{2\sqrt\gamma}e^{-|\cdot|/\sqrt\gamma}*f.
\]
For $C\geq2$ and $j=-1,0$,
\begin{align}
\|\chi_C R_\gamma f\|_{H^{j+2}}
&\leq C_0\gamma^{-1}\|\chi_C f\|_{H^j},\label{E:R-sobolev}\\
\|\chi_C R_\gamma\partial_yF\|_2
&\leq C_0\gamma^{-3/4}\|\chi_CF\|_1.\label{E:R-L1}
\end{align}
In addition,
\begin{equation}\label{E:R-small}
\|(R_\gamma-1)f\|_2\leq\frac{\sqrt\gamma}{2}\|f'\|_2,
\qquad
\|\sqrt\gamma\partial_yR_\gamma f\|_2\leq\frac12\|f\|_2.
\end{equation}
Their weighted versions hold with a fixed constant in place of $1/2$;
for the first, one can use $C_0\sqrt\gamma\|\chi_Cf\|_{H^1}$.
For a smooth bounded coefficient $b$,
\begin{equation}\label{E:R-comm}
[R_\gamma,b]R_\gamma^{-1}
=\gamma R_\gamma(2b'\partial_y+b'').
\end{equation}
\end{lemma}
\begin{proof}
The Fourier symbol $(1+\gamma k^2)^{-1}$ proves the unweighted
Sobolev estimates. The inequalities
$\gamma k^2/(1+\gamma k^2)\leq\sqrt\gamma|k|/2$ and
$\sqrt\gamma|k|/(1+\gamma k^2)\leq1/2$ give
\eqref{E:R-small}. The derivative of the convolution kernel has
$L^2$ norm $C\gamma^{-3/4}$, which proves \eqref{E:R-L1} by
Young's inequality. Weights can be inserted because
\[
\frac{\chi_C(x)}{\chi_C(y)}\leq e^{|x-y|/C},\qquad
|\partial_y^j\chi_C|\leq C_jC^{-j}\chi_C.
\]
The same kernel estimates apply to inverse weights. Duality and
commutation with derivatives give the estimate for $j=-1$.
For the weighted smallness bound, write
$R_\gamma-1=\gamma\partial_yR_\gamma\partial_y$ and move the
weight through the outer derivative. Finally the resolvent identity
$[R_\gamma,b]=\gamma R_\gamma[\partial_y^2,b]R_\gamma$
gives \eqref{E:R-comm}.
\end{proof}

\subsection{Flattening before regularization}
Suppress $\beta$ on operators within computations when no ambiguity
can arise. Set
\begin{equation}\label{E:P-flat}
\begin{aligned}
\mathsf P_\beta&=
\begin{bmatrix}\sqrt{\zeta_\beta}&0\\0&\zeta_\beta^{-1/2}\end{bmatrix},
\qquad r=(\log\zeta_\beta)',\\
\widetilde L_\beta&=\mathsf P_\beta\mathcal L_\beta\mathsf P_\beta
=\begin{bmatrix}
-\frac14\partial_y^2+U&-\sin\beta/\zeta_\beta\\
-\sin\beta/\zeta_\beta&1
\end{bmatrix},
\qquad U=\zeta_\beta+\frac18r'+\frac1{16}r^2.
\end{aligned}
\end{equation}
These multipliers and their inverses are bounded on all the Sobolev
spaces used here, uniformly on $I$. To verify the upper-left entry,
apply $\sqrt\zeta(-\frac14\partial_y\zeta^{-1}\partial_y)
\sqrt\zeta$ to a test function. Its first derivative term cancels,
and its zeroth order term is $-r'/8+r^2/16$; adding the zeroth order
part of $\mathcal L$ gives $U$.

For frozen $\beta$, the change $z=\mathsf P^{-1}v$ in the formal
dual equation gives
\begin{equation}\label{E:G-flat}
z_t=\mathscr G_\beta z,\qquad
\mathscr G_\beta=-2\sigma\widetilde L_\beta
\left(\partial_y+\frac r2\begin{bmatrix}1&0\\0&-1\end{bmatrix}\right).
\end{equation}
Indeed $\mathsf P'=\mathsf P\operatorname{diag}(r/2,-r/2)$ and
$\mathsf P^{-1}\sigma\mathcal L\mathsf P=\sigma\widetilde L$.
Its entries are
\begin{equation}\label{E:G-entries}
\mathscr G_\beta=
\begin{bmatrix}
\displaystyle\frac{2\sin\beta}{\zeta_\beta}\partial_y
+\frac{\sin\beta}{\zeta_\beta}r&-2\partial_y+r\\[3pt]
\displaystyle\frac12\partial_y^3+\frac14\partial_y^2(r\,\cdot)
-2U\partial_y-Ur&
\displaystyle\frac{2\sin\beta}{\zeta_\beta}\partial_y
-\frac{\sin\beta}{\zeta_\beta}r
\end{bmatrix}.
\end{equation}
The third order coefficient is constant. This is why regularization
is performed after the change of variables: the highest order term
then commutes with $R_\gamma$.

Define the one dual variable used in the nonlinear proof by
\begin{equation}\label{E:w-definition}
w=R_\gamma\mathsf P_\beta^{-1}\sigma\mathcal L_\beta
\left(\chi_A\begin{bmatrix}\hat\theta\\\hat\xi\end{bmatrix}\right),
\qquad
\mathcal D_w(t)=D_B(w(t)).
\end{equation}
The first component before regularization is in $L^2$, and the second
is in $H^{-1}$. Lemma~\ref{L:resolvent} gives
\begin{equation}\label{E:w-bound}
\|w_1\|_{H^2}+\|w_2\|_{H^1}\leq C_\gamma\mathcal D_A^{1/2}.
\end{equation}
Thus $w$ is defined for every energy perturbation. Localizing first
is important: the global norms of $w$ are controlled by
$\mathcal D_A$, which will appear with a favorable sign in the
primal virial estimate.

\begin{lemma}[Transferred constraint and its defect]\label{L:transferred}
The variable \eqref{E:w-definition} satisfies
\begin{equation}\label{E:w-constraint}
\langle w,R_\gamma^{-1}\mathsf P_\beta h_\beta\rangle
=\cos\beta\left\langle
\chi_A\begin{bmatrix}\hat\theta\\\hat\xi\end{bmatrix},
\sigma\begin{bmatrix}\hat\eta_\beta\\\hat\nu_\beta\end{bmatrix}\right\rangle.
\end{equation}
For either test vector $g=\sigma(\hat\eta_\beta,\hat\nu_\beta)^T$ or
$g=(\hat\eta_\beta',\hat\nu_\beta')^T$, the orthogonality conditions imply
\begin{equation}\label{E:constraint-defects}
\left|\left\langle
\chi_A\begin{bmatrix}\hat\theta\\\hat\xi\end{bmatrix},g\right\rangle\right|
\leq CA^{-2}\mathcal D_A^{1/2}.
\end{equation}
\end{lemma}
\begin{proof}
Move $R_\gamma$ to the smooth test, use
$\mathsf P^{-1}\sigma=\sigma\mathsf P$, and then use
\eqref{E:L-identities}. This gives \eqref{E:w-constraint} exactly.
For the second assertion, the original pairing with $g$ is zero, so
\[
\left\langle\chi_A\begin{bmatrix}\hat\theta\\\hat\xi\end{bmatrix},g\right\rangle
=\left\langle\chi_A\begin{bmatrix}\hat\theta\\\hat\xi\end{bmatrix},
(1-\chi_A^{-1})g\right\rangle.
\]
The exponentially decaying test satisfies
$\|(1-\cosh(y/A))g\|_2\leq CA^{-2}$.
Cauchy--Schwarz proves \eqref{E:constraint-defects}.
\end{proof}

\begin{lemma}[Recovery of the density and translation control]
\label{L:density-inverse}
For $B$ sufficiently large, $A\geq B$, and $0<\gamma\leq1$,
\begin{equation}\label{E:density-inverse}
\|\rho\hat\theta\|_2^2\leq C\mathcal D_w+CA^{-4}\mathcal D_A.
\end{equation}
Moreover,
\begin{equation}\label{E:alpha-w}
\begin{aligned}
|\dot a-2\sin\beta|
&\leq C\mathcal D_w^{1/2}+CA^{-2}\mathcal D_A^{1/2}+C\mathcal D_A,\\
|\dot\beta|&\leq C\mathcal D_A.
\end{aligned}
\end{equation}
These constants are independent of $A$ large and $\gamma\leq1$.
\end{lemma}
\begin{proof}
Put $v_A=\mathsf P(1-\gamma\partial_y^2)w$ in this proof. By
\eqref{E:w-definition} it equals
$\sigma\mathcal L(\chi_A\hat\theta,\chi_A\hat\xi)^T$.
Eliminating the second component gives
\begin{align*}
\chi_A\hat\xi&=\sin\beta\,\zeta_\beta^{-2}\chi_A\hat\theta
+\zeta_\beta^{-1}(v_A)_1,\\
H_\beta(\chi_A\hat\theta)
&=(v_A)_2+\sin\beta\,\zeta_\beta^{-2}(v_A)_1.
\end{align*}
The kernel of $H_\beta$ is $\spn\{\hat\eta_\beta'\}$.
With $g_\beta=\hat\eta_\beta'
+\sin\beta\,\zeta_\beta^{-2}\hat\nu_\beta'$ we have
$\langle\hat\eta_\beta',g_\beta\rangle
=\|(\hat\eta_\beta',\hat\nu_\beta')^T\|_2^2>0$.
Consequently the rank-one correction
\[
H_\beta+\hat\eta_\beta'\langle\,\cdot\,,g_\beta\rangle:
L^2\longrightarrow H^{-2}
\]
is an isomorphism. To check injectivity, pair its kernel equation
with $\hat\eta_\beta'$, then use the kernel of $H_\beta$.
Surjectivity follows by first solving on the complement of that
kernel and then fixing its scalar coefficient.
Conjugating by $\chi_B$ changes this operator by
$O(B^{-1})$ in $L^2\to H^{-2}$ norm. For the rank-one part this
uses the exponential decay of both test functions. A Neumann series
therefore gives the same inverse with weights when $B$ is large.

The remaining scalar coefficient is determined by translation:
\[
\langle\chi_A\hat\theta,g_\beta\rangle
=\left\langle\chi_A\begin{bmatrix}\hat\theta\\\hat\xi\end{bmatrix},
\begin{bmatrix}\hat\eta_\beta'\\\hat\nu_\beta'\end{bmatrix}\right\rangle
-\langle\zeta_\beta^{-1}(v_A)_1,\hat\nu_\beta'\rangle.
\]
The operator $1-\gamma\partial_y^2$, followed by multiplication by
the profile coefficients, maps weighted $L^2$ into weighted $H^{-2}$
uniformly for $\gamma\leq1$. The inverse just constructed and
\eqref{E:constraint-defects} thus yield
\[
\|\chi_B\chi_A\hat\theta\|_2
\leq C\|\chi_Bw\|_2+CA^{-2}\mathcal D_A^{1/2}.
\]
Taking $2/B<k$ and $A\geq B$ bounds $\rho/(\chi_A\chi_B)$ uniformly.
This proves
\eqref{E:density-inverse}.

For the translation estimate, move the linear functional in
\eqref{E:modulation-bounds} to
$\mathcal L(\chi_A\hat\theta,\chi_A\hat\xi)^T$, using the test
$\sigma(\hat\eta_\beta'',\hat\nu_\beta'')^T/\chi_A$.
The commutator $[\mathcal L,\chi_A]$ paired with this localized test
is $O(A^{-2}\mathcal D_A^{1/2})$; on the soliton core
$\chi_A'=O(|y|A^{-2})$. The main term becomes
\[
\left\langle w,
R_\gamma^{-1}\left(\frac{\mathsf P}{\chi_A}
\begin{bmatrix}\hat\eta_\beta''\\\hat\nu_\beta''\end{bmatrix}\right)
\right\rangle,
\]
which is bounded by $C\|\chi_Bw\|_2$ uniformly for $\gamma\leq1$.
The quadratic remainder in \eqref{E:modulation-bounds} is bounded
by $C\mathcal D_A$. This proves both estimates in \eqref{E:alpha-w}.
\end{proof}

\section{The primal localized virial}
\label{S:primal}
Define the bounded functional
\begin{equation}\label{E:primal-functional}
\mathcal I_A(t)=\int\zeta_{\beta(t)}\phi_A\hat\theta\hat\xi\,dy.
\end{equation}
The broad truncation guarantees $|\mathcal I_A|\leq C A\delta^2$
under orbital stability.

\begin{lemma}[Primal linear form]\label{L:primal-form}
Let $p=\zeta_\beta ^{-1}$, $M=-\sin\beta /\zeta_\beta $ and
$\mathcal B=1+\tfrac14\zeta_\beta ^{-2}\zeta_\beta ''
-\tfrac14\zeta_\beta ^{-3}(\zeta_\beta ')^2$. For a real smooth vector $F$, the
linear derivative of \eqref{E:primal-functional} is
\begin{align}\label{E:primal-form}
\mathfrak p_{\beta,A}(F)
={}&\int\frac{3p(\zeta_\beta\phi_A)'-p'\zeta_\beta\phi_A}{4}|F_1'|^2\,dy\\
&+\int\left(\mathcal B(\zeta_\beta\phi_A)'-\mathcal B'\zeta_\beta\phi_A
-\tfrac14(p(\zeta_\beta\phi_A)'')'\right)F_1^2\,dy\nonumber\\
&+2\int(M(\zeta_\beta\phi_A)'-M'\zeta_\beta\phi_A)F_1F_2\,dy
+\int(\zeta_\beta (\zeta_\beta\phi_A)'-\zeta_\beta '\zeta_\beta\phi_A)F_2^2\,dy.\nonumber
\end{align}
For $A$ sufficiently large,
\begin{equation}\label{E:primal-Garding}
\mathfrak p_{\beta,A}\!\left(\begin{bmatrix} \hat\theta \\ \hat\xi \end{bmatrix}\right)
\geq c_0\mathcal D_A-C\left\|\rho\hat\theta\right\|_2^2,
\end{equation}
with $c_0>0$ independent of $A$.
\end{lemma}
\begin{proof}
For $\partial_t (\hat\theta,\hat\xi)^T=J\mathcal L(\hat\theta,\hat\xi)^T$, integration by parts gives
$\mathcal I_A'=2\langle\mathcal L(\hat\theta,\hat\xi)^T,
\partial_y(\zeta_\beta\phi_A(\hat\theta,\hat\xi)^T)\rangle$ and \eqref{E:primal-form}.
The weight $\zeta_\beta\phi_A$ gives exactly
\begin{align*}
\frac{3p(\zeta_\beta\phi_A)'-p'\zeta_\beta\phi_A}{4}
&=\tfrac34\chi_A^2+\phi_A r,\\
\zeta_\beta(\zeta_\beta\phi_A)'-\zeta_\beta'\zeta_\beta\phi_A
&=\zeta_\beta^2\chi_A^2,\\
M(\zeta_\beta\phi_A)'-M'\zeta_\beta\phi_A
&=-\sin\beta(\chi_A^2+2\phi_A r).
\end{align*}
In particular the derivative coefficient and the $F_2^2$ coefficient
are positive. The normalized zero-order matrix tends to
$\left(\begin{smallmatrix}1&-\sin\beta \\-\sin\beta &1\end{smallmatrix}\right)$
outside the core, with an $O(A^{-2})$ cutoff error in its first entry.
More explicitly, reserve the fraction $\omega=\cos^{2}\beta /2$ of the positive
$F_2^2$ term and complete the square in the remainder. The resulting
$F_1^2$ coefficient at infinity is
$1-\sin^{2}\beta /(1-\omega)=\cos^{2}\beta /(1+\sin^{2}\beta )>0$.
Every other profile term is bounded by $C\rho^2F_1^2$;
this uses $|\phi_A|\leq|y|$, the derivative bounds for $\chi_A$, and
the exponential profile decay. Absorb the $O(A^{-2})$ term by choosing
$A$ large. Since $\phi_A r\geq0$, \eqref{E:primal-Garding} follows.
\end{proof}

\begin{lemma}[Primal nonlinear estimate]\label{L:primal-nonlinear}
For smooth small perturbations with $1-\eta\geq c_I>0$,
\begin{equation}\label{E:primal-nonlinear}
\left|2\left\langle\mathcal N_\beta,
\partial_y(\zeta_\beta\phi_A\begin{bmatrix} \hat\theta \\ \hat\xi \end{bmatrix})\right\rangle\right|
\leq C\delta \mathcal D_A.
\end{equation}
The constant is independent of sufficiently large $A$.
\end{lemma}
\begin{proof}
This is a variational cancellation, which can be made explicit using
\eqref{E:cubic-density}. In this proof only, let $m=\zeta_\beta\phi_A$
and regard $\mathcal R_\beta$ as a function of the four independent
arguments $(y,\hat\theta,\hat\theta_y,\hat\xi)$.
Its partial $y$ derivative differentiates the profile coefficients
only. Since $\mathcal N_\beta$ is its variational derivative, the
chain rule and integration by parts give
\begin{align*}
\left\langle\mathcal N_\beta,
\partial_y\left(m\begin{bmatrix}\hat\theta\\\hat\xi\end{bmatrix}\right)
\right\rangle
={}&\int m'\left(
\hat\theta\,\partial_{\hat\theta}\mathcal R_\beta
+2\hat\theta_y\,\partial_{\hat\theta_y}\mathcal R_\beta
+\hat\xi\,\partial_{\hat\xi}\mathcal R_\beta-\mathcal R_\beta\right)\\
&+\int m''\hat\theta\,\partial_{\hat\theta_y}\mathcal R_\beta
-\int m\,\partial_y^{\mathrm{explicit}}\mathcal R_\beta.
\end{align*}
Indeed the terms with coefficient $m$ combine into
$m(d_y\mathcal R_\beta-\partial_y^{\mathrm{explicit}}\mathcal R_\beta)$.
This identity removes both $\hat\xi_y$ and $\hat\theta_{yy}$.
The explicit expression \eqref{E:cubic-density}, uniform
nonvanishing, and $\|\hat\theta\|_\infty\leq C\delta$ show that
the first two integrands, without $m'$ and $m''$, are bounded by
$C\delta(\hat\theta_y^2+\hat\theta^2+\hat\xi^2)$.
The explicit $y$ derivative has an exponentially decaying coefficient
and the same bound. Finally,
\[
|m'|+|m''|\leq C\chi_A^2,\qquad
|m|\,|\partial_y^j\begin{bmatrix}\hat\eta_\beta\\\hat\nu_\beta\end{bmatrix}|
\leq C_j\chi_A^2\quad(j\geq1),
\]
uniformly for large $A$, since $|\phi_A|\leq|y|$.
The asserted estimate follows. Only the energy norm is used;
no $L^\infty$ bound on $\hat\theta_y$ or $\hat\xi$ enters.
\end{proof}

\begin{proposition}[Primal differential inequality]\label{P:primal}
For $B,A$ large, $A\geq B$, and $\delta$ sufficiently small,
\begin{equation}\label{E:primal-inequality}
\mathcal I_A'\geq c_1\mathcal D_A-C_1\mathcal D_w.
\end{equation}
$c_1>0$ and $C_1$ may be taken independent of $A$ and of $\gamma\leq1$.
\end{proposition}
\begin{proof}
Differentiate using the exact perturbation equation. Lemmas
\ref{L:primal-form}, \ref{L:primal-nonlinear}, and
\ref{L:density-inverse} control the linear and nonlinear terms.
The advection term is
$-(\dot a-2\sin\beta)\int(\zeta_\beta\phi_A)'\hat\theta\hat\xi=O(\delta \mathcal D_A)$.
The profile translation term is
$(\dot a-2\sin\beta)\langle\zeta_\beta\phi_A\sigma(\hat\eta_\beta',\hat\nu_\beta')^T,(\hat\theta,\hat\xi)^T\rangle$.
By \eqref{E:alpha-w}, it is bounded, for any $\varepsilon>0$, by
\[
\varepsilon \mathcal D_A+C_\varepsilon \mathcal D_w
+C A^{-2}\mathcal D_A+C\delta \mathcal D_A.
\]
The remaining terms come from the derivative of the soliton with
respect to $\beta$ and the time derivative of the weight. Their
pairings have size $C\delta$ and $C\delta^2$, respectively.
Since $|\dot\beta|\leq C\mathcal D_A$, both are
$O(\delta\mathcal D_A)$.
Choose $\varepsilon$ small, then $A$ large and $\delta$ small, and absorb
the resulting perturbations of $\mathcal D_A$.
\end{proof}

\section{The regularized dual localized virial}
\label{S:dual}
For $w$ in \eqref{E:w-definition}, define
\begin{equation}\label{E:dual-functional}
\mathcal V_{B,\gamma}(t)=\int\zeta_{\beta(t)}\phi_Bw_1w_2\,dy.
\end{equation}
Under orbital stability, $|\mathcal V_{B,\gamma}|\leq C_{B,\gamma}\delta^2$.

\begin{lemma}[Truncated quadratic coercivity]\label{L:truncated-coercivity}
Choose $B$ sufficiently large and then
$\gamma$ sufficiently small. For $A\geq B$ large,
\begin{equation}\label{E:truncated-coercivity}
\left\langle\mathcal W_{\beta,B}\mathsf P w,\mathsf P w\right\rangle
\geq\kappa_0\mathcal D_w-C A^{-4}\mathcal D_A,
\qquad
\mathcal W_{\beta,B}=\partial_y\mathcal L\zeta_\beta\phi_B
-\zeta_\beta\phi_B\mathcal L\partial_y,
\end{equation}
with $\kappa_0>0$ uniform on $I$.
\end{lemma}
\begin{proof}
The direct commutator expansion with the weight $\zeta_\beta\phi_B$ gives
\begin{align*}
(\mathcal W_{\beta,B})_{11}
&=-\partial_y\left(\tfrac34\chi_B^2+\tfrac12\phi_Br\right)\partial_y
+\zeta_\beta \chi_B^2+\zeta_\beta '\phi_B
-\tfrac14(\phi_B'''+\phi_B''r+\chi_B^2r'),\\
(\mathcal W_{\beta,B})_{12}&=-\sin\beta \chi_B^2,\qquad
(\mathcal W_{\beta,B})_{22}=\zeta_\beta ^2\chi_B^2+2\zeta_\beta \zeta_\beta '\phi_B.
\end{align*}
Compare this form with the untruncated form evaluated on $\chi_Bv$.
The inequalities
\[
|\phi_B'''|\leq C B^{-2}\chi_B^2,
\qquad |(\phi_B-y\chi_B^2)\partial_y^jr|
\leq C_j B^{-2}\chi_B^2
\]
and integration of the cross derivative show
\begin{equation}\label{E:weight-comparison}
\left|\left\langle\mathcal W_{\beta,B}v,v\right\rangle
-\left\langle\mathcal W_\beta(\chi_Bv),\chi_Bv\right\rangle\right|
\leq C B^{-2}D_B(v).
\end{equation}
For the second inequality, use the Taylor expansion on $|y|\leq B$
and the exponential decay of $r$ on its complement.

Theorem~\ref{T:virial-coercivity} also gives the penalty version
\[
\left\langle\mathcal W_\beta F,F\right\rangle
\geq\tfrac12\kappa_I\left\|F\right\|_{H^1\times L^2}^2
-C|\left\langle F,h_\beta\right\rangle|^2.
\]
It follows by projecting onto $h_\beta^\perp$ and bounding the cross
pairing of the fixed smooth profile by the form norm.
Take $F=\chi_B\mathsf P w$. Equations
\eqref{E:w-constraint}--\eqref{E:constraint-defects} imply
\[
|\left\langle F,h_\beta\right\rangle|
\leq C(\gamma+B^{-2})\mathcal D_w^{1/2}
+C A^{-2}\mathcal D_A^{1/2}.
\]
Indeed $R_\gamma^{-1}(\mathsf P h_\beta)=(\mathsf P h_\beta)-\gamma(\mathsf P h_\beta)''$,
and $(1-\chi_B^{-1})(\mathsf P h_\beta)$ is $O(B^{-2})$ in $L^2$.
Use boundedness of $\mathsf P^{\pm1}$ on the weighted form space,
then choose $B$ and $\gamma$ so that the $\mathcal D_w$ errors are absorbed.
\end{proof}

\begin{lemma}[Small linear regularization error]\label{L:flat-commutator}
For $B$ sufficiently large and $0<\gamma\leq\gamma_0(I,B)$,
\begin{equation}\label{E:flat-commutator}
\left|\left\langle\zeta_\beta\phi_B\sigma
[R_\gamma,\mathscr G_\beta]R_\gamma^{-1}w,w\right\rangle\right|
\leq C_B\sqrt\gamma \mathcal D_w.
\end{equation}
\end{lemma}
\begin{proof}
Write $m=\zeta_\beta\phi_B$, $X_j(f)=\|\chi_B f\|_{H^j}$ for
$j=0,1$, and $D_b=2b'\partial_y+b''$.  The constant-coefficient
third-order term commutes with $R_\gamma$.  For each of the other
coefficients $b$ in \eqref{E:G-entries}, the derivatives $b',b'',b'''$
are bounded by $C_B\chi_B^2$ when $B$ is large.  Multiplication by
$m$ is bounded in both $X_0$ and $X_1$.  By \eqref{E:R-comm}, a term
$b\partial_y$ gives the commutator $\gamma R_\gamma D_b\partial_y$.
Put $q=R_\gamma(mg)$. Since $D_b^*=-D_b$, integration by parts gives
\[
 \langle R_\gamma D_b\partial_y f,mg\rangle
 =-\langle f',D_b q\rangle
 =\langle f,2b'q''+3b''q'+b'''q\rangle .
\]
The weighted resolvent bounds give
$\gamma(\|\chi_B q'\|_2+\|\chi_B q\|_2)
 \le C_B\sqrt\gamma X_0(g)$ and, since
$\gamma\partial_y^2R_\gamma=R_\gamma-1$,
$\gamma\|\chi_B q''\|_2\le C_B\sqrt\gamma X_1(g)$.
Using $|b^{(j)}|\le C_B\chi_B^2$ to put one weight on each factor,
the first pairing is thus bounded, after multiplication by $\gamma$,
by $C_B\sqrt\gamma X_1(f)X_0(g)$ and the second by
$C_B\sqrt\gamma X_0(f)X_1(g)$. These two estimates apply to the
diagonal entries with $(f,g)=(w_1,w_2)$ and $(w_2,w_1)$,
respectively.  The first estimate also treats $-2U\partial_y$ in
the lower-left entry, where $f=g=w_1$.  A zeroth-order term $b$
gives $\gamma R_\gamma D_b$, whose pairing equals
$-\gamma\langle f,D_b q\rangle$ and is bounded by
$C_B\sqrt\gamma X_0(f)X_0(g)$.

The only remaining highest-order contribution, in the lower-left entry,
is exactly
\[
\tfrac14[R_\gamma,\partial_y^2(r\,\cdot)]R_\gamma^{-1}
=\tfrac14(R_\gamma-1)D_r.
\]
It pairs $w_1$ with itself. Put $T_\gamma=R_\gamma-1$.
For a real test function $f$, self-adjointness gives
\[
\left\langle\zeta_\beta\phi_B T_\gamma D_r f,f\right\rangle
=\left\langle D_r f,T_\gamma(\zeta_\beta\phi_Bf)\right\rangle.
\]
The weighted version of \eqref{E:R-small} implies
$\|\chi_BT_\gamma(\zeta_\beta\phi_Bf)\|_2
\leq C_B\sqrt\gamma\|\chi_Bf\|_{H^1}$.
Since $r'$ and $r''$ decay faster than $\chi_B^2$,
$\|\chi_B^{-1}D_rf\|_2\leq C_B\|\chi_Bf\|_{H^1}$.
Cauchy--Schwarz therefore bounds this contribution by
$C_B\sqrt\gamma\|\chi_Bf\|_{H^1}^2$.
All derivatives of $2\sin\beta/\zeta_\beta$, $U$, and $r$ in these commutators decay exponentially;
inserting $\chi_B^{-1}$ in the kernel estimates is therefore harmless
for $B$ large. The bounded logarithmic derivatives of $\chi_B$
produce only terms of the same or lower order. This proves
\eqref{E:flat-commutator} for the complete matrix.
\end{proof}

\begin{lemma}[Nonlinear source at energy regularity]\label{L:dual-source}
The regularized dual source $\mathfrak f$ arising from
$\chi_AJ\mathcal N_\beta$ satisfies
\begin{equation}\label{E:dual-source}
\left\|\mathfrak f_1\right\|_2+\left\|\mathfrak f_2\right\|_{H^{-1}}
\leq C_\gamma\delta \mathcal D_A^{1/2}.
\end{equation}
The source arising from $[\chi_A,J\mathcal L](\hat\theta,\hat\xi)^T$ satisfies
the same estimate with $\delta$ replaced by $A^{-1}$.
These sources mean applying
$R_\gamma\mathsf P^{-1}\sigma\mathcal L$ to the indicated
localized perturbation equation terms.
\end{lemma}
\begin{proof}
From \eqref{E:remainder-GF},
\[
\left\|\chi_A G\right\|_2+\left\|\chi_A F\right\|_1
+\left\|\chi_A\hat\theta\hat\xi\right\|_2
\leq C\delta \mathcal D_A^{1/2}.
\]
Here $\hat\theta\hat\theta'$ and $\hat\theta\hat\xi$ use
$\|\hat\theta\|_\infty\leq C\delta$, while, for example,
$\|\chi_A(\hat\theta')^2\|_1
\leq\|\hat\theta'\|_2\|\chi_A\hat\theta'\|_2$;
the same estimate applies to $\hat\xi^2$.
Write $\chi_AG'=\partial_y(\chi_AG)-\chi_A'G$.
The first source component contains at most two derivatives of
$\chi_AG$ and one derivative of $\chi_AF$.
The second contains, additionally, at most three derivatives of
$\chi_A\hat\theta\hat\xi$. Resolvent estimates give, respectively,
$L^2$ and $H^{-1}$ bounds; variable smooth profile coefficients
and derivatives of $\chi_A$ give lower-order terms of the same size.

For the linear localization source, use
\[
[\chi_A,J\mathcal L]\begin{bmatrix} \hat\theta \\ \hat\xi \end{bmatrix}
=2\sigma\bigl(\chi_A'\mathcal L\begin{bmatrix} \hat\theta \\ \hat\xi \end{bmatrix}
+\partial_y[\mathcal L,\chi_A]\begin{bmatrix} \hat\theta \\ \hat\xi \end{bmatrix}\bigr).
\]
The commutator $[\mathcal L,\chi_A]$ differentiates only $\hat\theta$ and
has at most one derivative; every coefficient has a factor
$O(A^{-1}\chi_A)$ or $O(A^{-2}\chi_A)$.
After applying the dual transform and $R_\gamma$, the largest terms
are resolvent-regularized second derivatives of $\chi_A'\hat\xi$ and
of a weighted $\hat\theta$. Their $L^2$ or $H^{-1}$ bounds are
$C_\gamma A^{-1}\mathcal D_A^{1/2}$.
This proves both assertions without requiring an unregularized dual
variable to belong to $H^1\times L^2$.
\end{proof}

\begin{proposition}[Dual differential inequality]\label{P:dual}
After choosing $B$ large and $\gamma$ small,
there is $\kappa_1>0$ such that, for $A\geq B$ large,
\begin{equation}\label{E:dual-inequality}
\mathcal V_{B,\gamma}'
\geq\kappa_1\mathcal D_w-C_{B,\gamma}(A^{-1}+\delta)\mathcal D_A.
\end{equation}
\end{proposition}
\begin{proof}
We first record the complete equation, to specify where each error
comes from. Within this proof write
$\mathcal T_\beta=\mathsf P_\beta^{-1}\sigma\mathcal L_\beta$;
this is an operator, not an additional perturbation variable.
Its coefficient derivatives are denoted by $\partial_y\mathcal T_\beta$
and $\partial_\beta\mathcal T_\beta$. Differentiating
\eqref{E:w-definition} and using \eqref{E:pert-eqn} gives
\begin{equation}\label{E:w-exact}
\begin{aligned}
w_t={}&\mathscr G_\beta w
+[R_\gamma,\mathscr G_\beta]R_\gamma^{-1}w\\
&+R_\gamma\mathcal T_\beta[\chi_A,J\mathcal L_\beta]
\begin{bmatrix}\hat\theta\\\hat\xi\end{bmatrix}
+R_\gamma\mathcal T_\beta(\chi_AJ\mathcal N_\beta)\\
&+(\dot a-2\sin\beta)w_y\\
&-(\dot a-2\sin\beta)R_\gamma
\left((\partial_y\mathcal T_\beta)\chi_A
+\mathcal T_\beta\chi_A'\right)
\begin{bmatrix}\hat\theta\\\hat\xi\end{bmatrix}\\
&+(\dot a-2\sin\beta)R_\gamma\mathcal T_\beta
\left(\chi_A\begin{bmatrix}\hat\eta_\beta'\\\hat\nu_\beta'\end{bmatrix}\right)\\
&-\dot\beta R_\gamma\mathcal T_\beta
\left(\chi_A\partial_\beta
\begin{bmatrix}\hat\eta_\beta\\\hat\nu_\beta\end{bmatrix}\right)
+\dot\beta R_\gamma(\partial_\beta\mathcal T_\beta)
\left(\chi_A\begin{bmatrix}\hat\theta\\\hat\xi\end{bmatrix}\right).
\end{aligned}
\end{equation}
The identities used here are
$\mathcal T_\beta J\mathcal L_\beta=\mathscr G_\beta\mathcal T_\beta$
and
$\chi_A\partial_yF=\partial_y(\chi_AF)-\chi_A'F$.
Thus no parameter derivative or localization term is omitted.

The homogeneous contribution to the derivative of
\eqref{E:dual-functional} is
$\langle\mathcal W_{\beta,B}\mathsf Pw,\mathsf Pw\rangle$.
Indeed, at fixed $\beta$, $\mathsf Pw$ obeys the formal dual equation,
and $(\mathsf Pw)_1(\mathsf Pw)_2=w_1w_2$.
Lemmas~\ref{L:truncated-coercivity} and \ref{L:flat-commutator}
therefore give a positive multiple of $\mathcal D_w$, with a loss
$CA^{-4}\mathcal D_A$, once $B$ is large and $\gamma$ small.

If a source $\mathfrak f$ is in $L^2\times H^{-1}$, its virial pairing is with
$(\zeta_\beta\phi_Bw_2,\zeta_\beta\phi_Bw_1)^T$ and has magnitude
at most
\[
C_B(\|\mathfrak f_1\|_2+\|\mathfrak f_2\|_{H^{-1}})\|w\|_{\mathcal X}.
\]
Consequently the two sources on the second line of
\eqref{E:w-exact} contribute at most
$C_{B,\gamma}(A^{-1}+\delta)\mathcal D_A$, by
Lemma~\ref{L:dual-source} and \eqref{E:w-bound}.
This estimate explains the broad localization before regularization:
both factors are bounded by the same $\mathcal D_A^{1/2}$.

Advection contributes exactly
\[
-(\dot a-2\sin\beta)\int(\zeta_\beta\phi_B)'w_1w_2,
\]
which is $O(\delta\mathcal D_w)$.
The coefficient-derivative and cutoff terms on the fourth line of
\eqref{E:w-exact} have $L^2\times H^{-1}$ norm at most
$C_\gamma|\dot a-2\sin\beta|\mathcal D_A^{1/2}$.
Here $\partial_y\mathcal T_\beta$ is still at most second order
and acts on a perturbation with one derivative in its first component.
They therefore contribute $O(C_{B,\gamma}\delta\mathcal D_A)$.

The profile translation term is smaller than its apparent size.
Without $\chi_A$ it is zero by \eqref{E:L-identities}. With $\chi_A$,
it is a smooth localized vector of norm at most
$C_\gamma A^{-2}|\dot a-2\sin\beta|$.
Pair it using \eqref{E:w-bound} and then \eqref{E:alpha-w}; for
any $\varepsilon>0$ its contribution is bounded by
\[
\varepsilon\mathcal D_w
+C_{B,\gamma,\varepsilon}A^{-4}\mathcal D_A
+C_{B,\gamma}\delta\mathcal D_A.
\]
The profile $\beta$ term is bounded by
$C_{B,\gamma}\delta|\dot\beta|
\leq C_{B,\gamma}\delta\mathcal D_A$.
The derivative of $\mathcal T_\beta$ has regularized norm at most
$C_\gamma|\dot\beta|\mathcal D_A^{1/2}$, so its pairing is
$O(C_{B,\gamma}\delta^2\mathcal D_A)$.
Finally differentiation of the weight contributes
\[
\dot\beta\int(\partial_\beta\zeta_\beta)\phi_Bw_1w_2
=O(C_{B,\gamma}\delta^2\mathcal D_A).
\]
Choose $\varepsilon$ small, then $A$ sufficiently large, and finally
$\delta$ small. Absorbing the indicated $\mathcal D_w$ errors proves
\eqref{E:dual-inequality}.
\end{proof}

\section{Closing the spacetime estimate}
\label{S:closure}

\begin{theorem}[Direct spacetime bound]\label{T:spacetime}
The choices
\eqref{E:choice-order} can be made so that every sufficiently small
energy-class perturbation obeys, for all $T>0$,
\begin{equation}\label{E:spacetime}
\sup_{0\leq t\leq T}\left\|\begin{bmatrix} \hat\theta(t) \\ \hat\xi(t) \end{bmatrix}\right\|_{\mathcal X}^2
+\int_0^T(\mathcal D_A+\mathcal D_w)\,dt\leq C\delta^2.
\end{equation}
The constant is independent of $T$. In addition,
\begin{equation}\label{E:modulation-integrability}
\int_0^\infty|\dot\beta|\,dt
+\int_0^\infty|\dot a-2\sin\beta|^2\,dt\leq C\delta^2.
\end{equation}
\end{theorem}
\begin{proof}
The orbital bound is Proposition~\ref{P:orbital}.
Choose $B$ large and $\gamma$ small so that Proposition~\ref{P:dual}
holds. In Proposition~\ref{P:primal}, write its constants as $c_1,C_1$.
Choose $M>0$ with $M\kappa_1-C_1\geq1$ and put
$\mathcal F=\mathcal I_A+M\mathcal V_{B,\gamma}$.
Next choose $A$ so large, and then $\delta$ so small, that
$M C_{B,\gamma}(A^{-1}+\delta)\leq c_1/2$.
Then
\[
\mathcal F'\geq\tfrac12c_1\mathcal D_A+\mathcal D_w,
\qquad |\mathcal F(t)|\leq C_{A,B,\gamma}\delta^2.
\]
Integrate from $0$ to $T$. The endpoint bound is independent of $T$,
which proves \eqref{E:spacetime}. Equation \eqref{E:alpha-w},
$\mathcal D_A\leq C\delta^2$, and the modulation equations then give
\eqref{E:modulation-integrability}.
Every error in this integration is proportional to one of the
controlled local energies. This is what makes the constant independent
of $T$ and permits the argument to be used on a prescribed finite
time interval.
\end{proof}

\subsection{Energy-class justification}
The differentiations above are first made for smooth initial data.
Approximate the nonvanishing hydrodynamic initial datum by smooth data
in $H^1\times L^2$, preserving a positive lower bound for $1-\eta$.
By continuous dependence for the hydrodynamic flow
\cite[Theorem~1.2]{Mohamad2014}, the approximating solutions converge
in $C([0,T];H^1\times L^2)$ on each finite interval. The orbital
bound gives a common positive lower bound on $1-\eta$ for all
sufficiently close approximants. In applying the smooth estimates
their initial distances can be taken to be $\delta+o(1)$, with the
smallness threshold chosen inside the orbital neighborhood.
The modulation maps are continuous and differentiable against their
smooth profile test functions. Their equations have a uniformly
invertible two-by-two coefficient matrix. Thus the parameters and
their weakly derived equations pass to the limit on finite intervals.

The decomposition $\mathcal N_1=G'+F$ is continuous from the energy
space into $H^{-1}$: $G$ converges in $L^2$, and $F$ in $L^1$.
For fixed $A,B,\gamma$, the regularized variable converges in
$C([0,T];H^1\times L^2)$ and both bounded virial functionals converge.
The integrated inequalities therefore pass to the energy solution;
alternatively their positive terms pass by lower semicontinuity.
All constants depend on the energy smallness and fixed cutoffs, not
on higher norms of the approximating data. Taking $T\to\infty$
proves the infinite-time statements. This use of approximation is
not a compactness or rigidity argument for an asymptotic limiting orbit.

\section{Asymptotic conclusions and proof of the main theorem}
\label{S:asymptotics}

\begin{lemma}[An elementary decay criterion]\label{L:decay-criterion}
If $g:[0,\infty)\to[0,\infty)$ is uniformly continuous and
$\int_0^\infty g(t)\,dt<\infty$, then $g(t)\to0$.
\end{lemma}
\begin{proof}
Otherwise uniform continuity gives disjoint intervals of a fixed
positive length on which $g$ is bounded below by a fixed positive
number, contradicting integrability.
\end{proof}

\begin{proposition}[Consequences of the spacetime estimate]\label{P:energy-asymptotics}
Under the hypotheses of Theorem~\ref{T:main}, there is
$\beta_\infty\in I$ with $\beta(t)\to\beta_\infty$ and
$|\beta_\infty-\beta(0)|\leq C\delta^2$.
For each $R>0$,
\begin{align}
\hat\theta(t)&\longrightarrow0\quad\text{in }L^2(-R,R),
\label{E:theta-local}\\
\hat\xi(t)&\rightharpoonup0\quad\text{in }L^2(\mathbb R),
\label{E:xi-weak}\\
\int_t^{t+1}\int_{-R}^R
(\hat\theta_y^2+\hat\theta^2+\hat\xi^2)(y,s)\,dy\,ds&\longrightarrow0.
\label{E:local-average}
\end{align}
There is also $b\in C^1([0,\infty))$ with
$b(t)-a(t)\to0$ and $b'(t)\to2\sin\beta_\infty$.
The statements \eqref{E:theta-local}--\eqref{E:local-average}
hold in the frame of $b$ as well.
\end{proposition}
\begin{proof}
The first assertion follows from \eqref{E:modulation-integrability}.
The exact hydrodynamic equation, the orbital bound, and the nonlinear
remainder estimates imply
\[
\sup_t\left\|\hat\theta_t(t)\right\|_{H^{-1}}<\infty,
\qquad \sup_t\left\|\hat\xi_t(t)\right\|_{H^{-2}}<\infty.
\]
For the second estimate use $\mathcal N_1=G'+F$ with
$F\in L^1\subset H^{-1}$ and $G\in L^2$.
For the first, the flux is a sum of $L^2$ terms, including $\hat\theta\hat\xi$.
Interpolation between $H^{-1}$ and the uniform $H^1$ bound makes
$\hat\theta$ uniformly continuous in $L^2$. Apply Lemma
\ref{L:decay-criterion} to $\|\chi_A\hat\theta\|_2^2$, whose integral
is finite by \eqref{E:spacetime}. This proves \eqref{E:theta-local}.
For every smooth compactly supported $\psi$, the scalar
$\langle\hat\xi(t),\psi\rangle$ is uniformly Lipschitz and belongs to
$L^2_t$. It therefore tends to zero. Density and the uniform $L^2$
bound give \eqref{E:xi-weak}. Equation \eqref{E:local-average}
follows directly from the tail of the spacetime integral.

We use the center smoothing of \cite[Lemma~2]{BGS_Asymp_Dark_2015}.
Put $v_\infty=2\sin\beta_\infty$. Choose a smooth nonnegative
$\omega$ supported in $(0,1)$ with integral one, and define
\[
b(t)=\int_0^1\omega(s)a(t+s)\,ds
-v_\infty\int_0^1s\omega(s)\,ds.
\]
Since $\dot a-2\sin\beta\in L^2_t$ and
$\beta(t)\to\beta_\infty$, Cauchy--Schwarz on $[t,t+1]$ gives
$b-a\to0$ and $b'\to v_\infty$.
The change of frame tends to zero: it does not affect local density
convergence or weak velocity convergence. For the time average, enlarge
the spatial interval by one for all large $t$ and use
\eqref{E:local-average} and smoothness of the profiles.
\end{proof}

\begin{proposition}[Strong local convergence at energy regularity]
\label{P:strong-asymptotics}
Under the assumptions of Theorem~\ref{T:main}, for every $R>0$,
\begin{equation}\label{E:strong-local}
 \|\hat\theta(t)\|_{H^1(-R,R)}+\|\hat\xi(t)\|_{L^2(-R,R)}\longrightarrow0,
 \qquad \dot a(t)-2\sin\beta(t)\longrightarrow0 .
\end{equation}
\end{proposition}
\begin{proof}
We first record a consequence of the conserved energies for an arbitrary
finite-energy solution $u$ of \eqref{E:GP}. In laboratory coordinates put
\[
 U(t)=(|u(t)|^2-1,u_x(t)),\qquad N=\|U(t)\|_2^2 .
\]
The number $N$ is twice the conserved physical energy. For a
nonvanishing solution it equals $2H_0(\eta(0),\nu(0))$ by
\eqref{E:conserved}. Using the unitary Fourier transform
componentwise, define, for $\tau\ge2$,
\[
 F_\tau(t)^2=\int_{\mathbb R}
       \frac{|\widehat U(t,\xi)|^2}{\tau^2+\xi^2}\,d\xi,\qquad
 D_\tau(t)=N-\tau^2 F_\tau(t)^2
     =\int_{\mathbb R}\frac{\xi^2}{\tau^2+\xi^2}
                    |\widehat U(t,\xi)|^2\,d\xi .
\]
Koch and Liao \cite[Theorem~1.6, (1.23)]{KochLiaoLow}, applied at $s=0$
in the normalization \eqref{E:GP}, give a conserved quantity
$\mathcal E^0_\tau$ satisfying
\[
 |\mathcal E^0_\tau(u(t))-F_\tau(t)^2|
       \le C\tau^{-1/2} F_\tau(t)^3
 \quad\text{when }\tau^{-1/2}F_\tau(t)\text{ is sufficiently small}.
\]
The finite-energy solution is also the solution in their $X^0$ flow
by consistency and uniqueness (see also \cite[Remark~3.3]{Wegner2023}).
Since $F_\tau(t)\le\sqrt{N}/\tau$, the smallness condition holds at
all times once $\tau$ is sufficiently large. Comparing the errors at
$0$ and $t$ gives
\[
 |D_\tau(t)-D_\tau(0)|\le 2C\tau^{-3/2}N^{3/2}.
\]
The multiplier in $D_\tau$ is at least $1/2$ for $|\xi|\ge\tau$,
and $D_\tau(0)\to0$ by dominated convergence. Hence
\begin{equation}\label{E:uniform-Fourier-tails}
 \tfrac12\sup_{t\ge0}\int_{|\xi|\ge\tau}|\widehat U(t,\xi)|^2\,d\xi
 \le D_\tau(0)+2C\tau^{-3/2}N^{3/2}\longrightarrow0 .
\end{equation}

Let $T_h f(y)=f(y+h)$. Plancherel's theorem, splitting at a fixed
frequency and using \eqref{E:uniform-Fourier-tails} on its complement,
gives $\sup_{t\ge0}\|T_h u_x(t)-u_x(t)\|_2\to0$ as $h\to0$.
Also $\|T_hu(t)-u(t)\|_\infty\le\sqrt N\,|h|^{1/2}$. In the gray
tube, $u$ and $1/u$ are uniformly bounded, so $1/u$ has the same
$O(|h|^{1/2})$ translation bound. It follows from
\[
 \eta_x=-2\Re(\bar u u_x),\qquad \nu=-\Im(u_x/u)
\]
that $\eta_x$ and $\nu$ are uniformly continuous under spatial
translations in $L^2$, uniformly in time. For example,
\[
 \|T_h(u_x/u)-u_x/u\|_2
 \le \|1/u\|_\infty\|T_h u_x-u_x\|_2
       +\|u_x\|_2\|T_h(1/u)-1/u\|_\infty .
\]
Translations by $a(t)$ preserve these bounds. The profile derivatives
$\hat\eta_\beta'$ and velocities $\hat\nu_\beta$, $\beta\in I$,
range over compact subsets of $L^2$. Subtracting them therefore
preserves uniform continuity under spatial translations.

After multiplication by a smooth compactly supported cutoff, the
families $\hat\theta_y(t)$ and $\hat\xi(t)$ are relatively compact
in $L^2$ by the Kolmogorov--Riesz theorem. Equation~\eqref{E:theta-local}
gives $\hat\theta_y(t)\to0$ in distributions, and
\eqref{E:xi-weak} identifies the weak velocity limit as zero.
Every sequence $t_n\to\infty$ thus has a subsequence along which both
localized quantities tend to zero strongly. Taking the cutoff equal to
one on $[-R,R]$ proves the local norm limits in \eqref{E:strong-local}.

Finally, the linear functional in \eqref{E:modulation-bounds} pairs
the perturbation with a uniformly exponentially localized profile test.
On $[-R,R]$ the local norm limits imply that the linear pairing
and the quadratic remainder both tend to zero. The global orbital
bound makes their tails uniformly small as $R\to\infty$. The denominator $\|(\hat\eta'_\beta,\hat\nu'_\beta)\|_2^2$
has a positive minimum on $I$. This gives the fitted-speed assertion.
\end{proof}

\begin{proof}[Proof of Theorem~\ref{T:main}]
Proposition~\ref{P:orbital} and Theorem~\ref{T:spacetime} give
\eqref{E:main-estimate}. Proposition~\ref{P:energy-asymptotics} gives
$\beta(t)\to\beta_\infty$ and the bound on $|\beta_\infty-\beta(0)|$.
Proposition~\ref{P:strong-asymptotics} now proves \eqref{E:main-strong}.
Smooth dependence of the profiles on $\beta\in I$ permits replacing
$\beta(t)$ by $\beta_\infty$ in the local limits. Local convergence
and the uniform bound in $\mathcal X$ imply weak convergence globally
by density of compactly supported tests.
\end{proof}

\begin{corollary}[Local convergence of the complex gradient in the
B\'ethuel--Gravejat--Smets frame]\label{C:BGS-gradient}
Under the hypotheses and in the notation of
\cite[Theorem~2]{BGS_Asymp_Dark_2015}, let $b(t)$ and $\theta(t)$
be the center and phase from that theorem and $U_{c^*}$ its limiting
soliton. For every $R>0$,
\[
 e^{-i\theta(t)}\Psi(\cdot+b(t),t)\longrightarrow U_{c^*}
      \quad\text{strongly in }H^1(-R,R)\quad(t\to\infty).
\]
\end{corollary}
\begin{proof}
Their theorem gives local uniform convergence of the fields and weak
$L^2(\mathbb R)$ convergence of the phase-corrected gradients to
$U'_{c^*}$. The high-frequency estimate
\eqref{E:uniform-Fourier-tails} applies to $\Psi_x$ after the fixed
rescaling between the two normalizations described in
Section~\ref{S:modulation}. Its Fourier tail is invariant under
translation and multiplication by a spatially constant phase, even
when the center and phase vary with time. Thus the corrected
gradients are relatively compact in $L^2(-R,R)$, by the cutoff
argument above. Their weak limit identifies every strong local limit.
This proves the assertion and answers the question following Theorem~2
of \cite{BGS_Asymp_Dark_2015}.
\end{proof}

The use of conserved energies here complements the earlier
integer-order identities of B\'ethuel, Gravejat, Saut and Smets
\cite{BGSSKdV2009} and the higher-regularity family of Koch and Liao
\cite{KochLiao2021}. For propagation of frequency equicontinuity by
conservation laws in the zero-background cubic NLS and mKdV equations,
see Harrop-Griffiths, Killip and Vi\c{s}an
\cite[Proposition~4.6]{HarropKillipVisan2024}.

\appendix
\section{Polynomial certificate for the virial form}\label{S:virial-certificate}
This appendix proves Lemma~\ref{L:virial-halfline}.
All notation introduced here is local to the calculation. Write $a,b,d$
for the coefficient functions $a_\beta,b_\beta,d_\beta$ in
\eqref{E:virial-coefficients}, and $r$ for the residual $r_\beta$ in
\eqref{E:virial-comparison-vector}.

\subsection{Reduction to rational functions}\label{S:virial-rational}
Here $x=y\cos\beta>0$, $t=\tanh x$, $v=t^2$, $e=1-v$, $p=\cos^2\!\beta$, and $Z=1-pe$. Subscripts $x,t$ indicate differentiation in these variables. The following formulas allow the calculation to be checked directly against \eqref{E:W-compact}:
\begin{align*}
 \rho&=(\log Z)_x=2pte/Z,&\rho_x&=2pe(1-3v)/Z-\rho^2,\\
 a&=3/4+x\rho/2,&a_x&=(\rho+x\rho_x)/2,\\
 b&=Z+2pxte-p\rho_x/4,&d&=Z(Z+4pxte).
\end{align*}
For $\ell=\log\psi_\beta$, we have $\ell=\log t-\frac13\log(1-p+3pv)+\frac{7}{16}\log e$. Consequently
\begin{equation}
 \mathcal H:=\Lambda_\beta/p=-a\big[e^2(\ell_{tt}+\ell_t^2)-2te\ell_t\big]-a_xe\ell_t+\frac{b-(1-p)/d}{p}.
 \label{E:certificate-H}
\end{equation}

To compute the residual, use $\partial_\beta\hat\eta_\beta=-2\sin\beta\cos\beta\,e(1-xt)$ and $\partial_\beta\hat\nu_\beta=\cos\beta\,\hat\eta_\beta/Z+\sin\beta\,\partial_\beta\hat\eta_\beta/Z^2$. Using \eqref{E:W-compact} in the residual definition \eqref{E:virial-comparison-vector} gives
\begin{align}
 g_1:=\frac{r_1}{\cos\beta\,e}&=\frac{p}{Z}-\frac{2(1-p)(1-xt)}{Z^2}
 -\frac{2}{p}\left[-ap(4-6e)+2pt a_x+b-\frac{1-p}{Z}\right],\nonumber\\
 g_2:=\frac{r_2}{\sin\beta\cos\beta}&=-2e(1-xt)-\frac{2}{p^2}\left[-pe+\frac{dpe}{Z}\right].\label{E:certificate-residual}
\end{align}
For example, $\partial_x\hat\eta_\beta=-2t\hat\eta_\beta$ and $\partial_x^2\hat\eta_\beta=(4-6e)\hat\eta_\beta$. In these variables \eqref{E:virial-pointwise} becomes
\begin{equation}
 \mathcal H\ge\tfrac14,\qquad G_1/\mathcal H+G_2\le15/2,\qquad
 G_1=e\left(g_1+\frac{(1-p)g_2}{de}\right)^2,\quad G_2=\frac{p(1-p)g_2^2}{de}.
 \label{E:certificate-target}
\end{equation}
Here $G_1,G_2$ denote scalar coefficients only in this calculation.

Put $z=e^{-x/2}$ and $j=pv/Z$. Then $t=(1-z^4)/(1+z^4)$, $p=j/(v+je)$, and $0<z,j<1$. The only nonrational dependence is $x=-2\log z$. The elementary bounds $4\tanh(x/4)\le x\le2\sinh(x/2)$ give
\[
 \frac{4(1-z)}{1+z}\le x\le z^{-1}-z.
\]
After multiplication by $z/t$, the ratio between the lower and upper endpoints is $4z/(1+z)^2$. Hence, for some $0\le\tau\le1$,
\begin{equation}
 \frac{zx}{t}=\frac{1+z^4}{1+z^2}
 \left[1-\tau\left(\frac{1-z}{1+z}\right)^2\right].\label{E:certificate-enclosure}
\end{equation}
We prove \eqref{E:certificate-target} on the entire cube $0\le z,j,\tau\le1$ after denominators are cleared. The substitution for $x$ is made \emph{after} the derivatives in \eqref{E:certificate-H} are taken.

\subsection{Exact polynomial certificates}\label{S:virial-polynomials}
Define the following positive factors on the interior of the cube:
\begin{align*}
 D&=(1+z^4)^2,&V&=(1-z^4)^2,&K&=(1+z^4)(1+z^2)(1+z)^2,\\
 N&=K(1+4je\,x/t),&&&D_*&=256KDVN(1+2j)^2.
\end{align*}
Using \eqref{E:certificate-H}--\eqref{E:certificate-enclosure}, clear denominators as follows:
\begin{align}
 J&=D_*\mathcal H,\qquad P_H=4J-D_*,\nonumber\\
 P_R&=D^2KN\big[(15-2G_2)J-2D_*G_1\big].\label{E:certificate-polynomials}
\end{align}
These are polynomials in $(z,j,\tau)$ with integer coefficients; a factored expansion is given below. Nonnegativity of $P_H$ gives $\mathcal H\ge1/4$, hence $J>0$ in the physical domain. Nonnegativity of $P_R$ then gives the second inequality in \eqref{E:certificate-target}.

For completeness, the two numerators in \eqref{E:certificate-polynomials} may be formed without rational-function arithmetic. Write
\begin{align*}
 B&=(1+z)^2-\tau(1-z)^2,& E&=4z^4,&E_2&=4z^2,\\
 M&=4z^3B,& N&=K+4jM,& A&=3K+4jM,\\
 b_0&=1+2j,& L&=D-3V-2jE,& P&=V+jE,\\
 Z_0&=zK+B(4E-D),& Y&=4MV-EK.&&
\end{align*}
Thus $e=E/D$, $v=V/D$, $ex/t=M/K$ and $a=A/(4K)$. Substitution and collection give
\begin{align*}
 C={}&(-2D+3jV+2j^2E-4jE+2E)K-A(D-3V)-4M(jL+V),\\
 R={}&zCN+2P\bigl[((1-j)DB-zK)N-(1-j)Z_0K\bigr],\\
 T={}&256KNb_0^2(P-E)+256Kb_0^2(1-j)Y+512Nb_0^2MV\\
     &\hspace{7mm}{}+AN\bigl[240Eb_0+b_0^2(56E-49V)\bigr],\\
 S={}&256Kb_0^2(1-j)EY-128KNb_0^2EL+128ANE^2(3-2j)\\
     &\hspace{7mm}{}-32Nb_0(EK+ML)(8E-7Vb_0),\\
 J={}&VT+jS,\qquad P_H=4J-256KDVNb_0^2,\\
 P_R={}&\bigl(15D^2KN-8E_2Vj(1-j)Z_0^2\bigr)J-512E_2Vb_0^2R^2.
\end{align*}
In particular $G_1=E_2R^2/(D^3K^2N^2)$ and $G_2=4E_2Vj(1-j)Z_0^2/(D^2KN)$. The signs can be read off in the tensor Bernstein basis. Explicitly, if $P=\sum_{\alpha}c_\alpha X^\alpha$ has multidegree $n$, its Bernstein coefficient at $i$ is
\begin{equation}
 B_i=\sum_{\alpha\le i} c_\alpha\prod_{k=1}^3\frac{\binom{i_k}{\alpha_k}}{\binom{n_k}{\alpha_k}}.\label{E:certificate-Bernstein}
\end{equation}
Here the tensor basis consists of
\[
\prod_{k=1}^3\binom{n_k}{i_k}X_k^{i_k}(1-X_k)^{n_k-i_k},
\qquad 0\leq X_k\leq1,
\]
and the polynomial is the sum of these basis functions with
coefficients $B_i$.
In each row below every coefficient is nonnegative. The last column is the floor of the least strictly positive coefficient; the preceding column counts zeros.
\begin{center}\footnotesize\setlength{\tabcolsep}{4pt}
\begin{tabular}{@{}llrrr@{}}\toprule
Polynomial & Multidegree & Count & Zeros & Least positive\\[-2pt]
 & & & & (floored)\\\midrule
$P_H(z,j,\tau)$ & $(32,5,2)$ & 594 & 6 & 180\\
$P_R(z/2,j,\tau)$ & $(64,8,4)$ & 2925 & 0 & 1635\\
$P_R((1+z)/2,j,\tau)$ & $(64,8,4)$ & 2925 & 10 & 849986\\\bottomrule
\end{tabular}\end{center}
The table is finite integer arithmetic: expand the displayed polynomials, substitute either $z/2$ or $(1+z)/2$ for the last two rows, and apply \eqref{E:certificate-Bernstein} in each variable. The companion script \path{verify_virial_certificate.py}, included with
the source files, also verifies \eqref{E:certificate-polynomials} directly from \eqref{E:certificate-H}--\eqref{E:certificate-residual} by formal rational-function arithmetic, then reproduces all three rows. Since the Bernstein basis functions are nonnegative on the closed cube, the table proves $P_H,P_R\ge0$. This proves Lemma~\ref{L:virial-halfline} and completes the proof.

\bibliographystyle{amsalpha}
\bibliography{references}

\end{document}